\documentclass[12pt]{amsproc}

\usepackage[latin1]{inputenc}
\usepackage{xspace,amssymb,amsfonts,euscript}
\usepackage{amsthm,amsmath}
\usepackage{palatino}
\usepackage{euscript}
\input xy \xyoption {all}

\RequirePackage{color}
\definecolor{myred}{rgb}{0.75,0,0}
\definecolor{mygreen}{rgb}{0,0.5,0}
\definecolor{myblue}{rgb}{0,0,0.65}

\RequirePackage{ifpdf}
\ifpdf
  \IfFileExists{pdfsync.sty}{\RequirePackage{pdfsync}}{}
  \RequirePackage[pdftex,
   colorlinks = true,
   urlcolor = myblue, 
   citecolor = mygreen, 
   linkcolor = myred, 
   pagebackref,
   bookmarksopen=true]{hyperref}
\else
  \RequirePackage[hypertex]{hyperref}
\fi

\RequirePackage{ae, aecompl, aeguill} 

    \def\AM{{\mathbb{A}}}
  \def\bg{{\mathfrak b}}  
    \def\CM{{\mathbb{C}}}
    \def\DM{{\mathbb{D}}}
    
    \def\FM{{\mathbb{F}}}
  \def\gg{{\mathfrak g}}

    \def\KM{{\mathbb{K}}}

    \def\NM{{\mathbb{N}}}
    \def\OM{{\mathbb{O}}}
\def\PG{{\mathfrak P}}    \def\PM{{\mathbb{P}}}
    \def\QM{{\mathbb{Q}}}
    \def\RM{{\mathbb{R}}}
\def\SG{{\mathfrak S}}    
  \def\tg{{\mathfrak t}}  
  \def\ug{{\mathfrak u}}

    \def\ZM{{\mathbb{Z}}}

    \def\BC{{\mathcal{B}}}
    \def\CC{{\mathcal{C}}}
    \def\DC{{\mathcal{D}}}
    
    \def\FC{{\mathcal{F}}}
    
    \def\HC{{\mathcal{H}}}
    \def\IC{{\mathcal{I}}}
    
    \def\KC{{\mathcal{K}}}
    \def\LC{{\mathcal{L}}}
    
    \def\NC{{\mathcal{N}}}
    \def\OC{{\mathcal{O}}}

  \def\sb{{\mathbf s}}

\def\FS{{\EuScript F}}

\def\KS{{\EuScript K}}

\def\SS{{\EuScript S}}

\def\a{\alpha}

\def\G{\Gamma}

\def\D{\Delta}
\def\e{\varepsilon}

\def\l{\lambda}

\newcommand{\nc}{\newcommand} \newcommand{\renc}{\renewcommand}

\newcommand{\rdots}{\mathinner{ \mkern1mu\raise1pt\hbox{.}
    \mkern2mu\raise4pt\hbox{.}
    \mkern2mu\raise7pt\vbox{\kern7pt\hbox{.}}\mkern1mu}}

\def\mini{{\mathrm{min}}}
\def\reg{{\mathrm{reg}}}
\def\subreg{{\mathrm{subreg}}}
\def\res{{\mathrm{res}}}
\def\rs{{\mathrm{rs}}}

\def\modules{{\text{-mod}}}

\DeclareMathOperator{\Ad}{Ad}

\DeclareMathOperator{\Tr}{Tr}
\DeclareMathOperator{\Irr}{Irr}
\DeclareMathOperator{\Ker}{Ker}
\DeclareMathOperator{\Coker}{Coker}

\def\pr{{\mathrm{pr}}}
\newcommand{\elem}[1]{\stackrel{#1}{\longto}}
\DeclareMathOperator{\RHOM}{R\underline{Hom}}

\def\ov{\overline}
\def\un{\underline}

\def\to{\rightarrow}

\def\longto{\longrightarrow}

\def\injto{\hookrightarrow}

\nc{\Br}{\mathcal{B}}
\nc{\id}{id}
\nc{\HotRR}{{}_R\mathcal{K}_R}
\nc{\HotR}{\mathcal{K}_R}
\nc{\excise}[1]{}
\nc{\defect}{\text{df}}
\nc{\h}[1]{\underline{H}_{#1}}

\nc{\Ga}{\mathbb{G}_a} 
\nc{\Gm}{\mathbb{G}_m} 

\nc{\Perv}{{\mathbf{P}}}

\nc{\IH}{{\mathrm{IH}}}

\nc{\ic}{\mathbf{IC}}

\nc{\gl}{{\mathfrak{gl}}}
\renc{\sl}{{\mathfrak{sl}}}
\renc{\sp}{{\mathfrak{sp}}}

\DeclareMathOperator{\codim}{{\mathrm{codim}}}
\DeclareMathOperator{\im}{{\mathrm{Im}}}
 \DeclareMathOperator{\Hom}{Hom}

 \DeclareMathOperator{\Rep}{Rep}

\newtheorem{thm}{Theorem}[section]

\newtheorem{prop}[thm]{Proposition}

\theoremstyle{definition}

\theoremstyle{remark}
\newtheorem{remark}[thm]{Remark}

\DeclareMathOperator{\Ext}{Ext}

\DeclareMathOperator{\Spec}{Spec}

\newcommand{\xra}{\xrightarrow}

\def\Iff{\Longleftrightarrow}

\def\uk{\underline{k}}

\def\pt{{\mathrm{pt}}}
\def\op{{\mathrm{op}}}

\DeclareMathOperator{\RGa}{R\G}
\DeclareMathOperator{\RGc}{R\G_c}
\DeclareMathOperator{\RHom}{RHom}
\DeclareMathOperator{\Gras}{Gr}

\def\Gr{{\EuScript Gr}}
 
\begin{document}

\title{Perverse sheaves and modular representation theory}

\author{Daniel Juteau} \address{ LMNO, Université de Caen
  Basse-Normandie, CNRS, BP 5186, 14032 Caen, France}
\email{juteau@math.unicaen.fr}

\author{Carl Mautner} \address{ Mathematics department, University of Texas at Austin, 1
  University Station C1200, Austin TX, 78712 USA}
\email{cmautner@math.utexas.edu}

\author{Geordie Williamson} \address{ Mathematical Institute,
  University of Oxford, 24-29 St Giles', Oxford, OX1 3LB, UK}
\email{williamsong@maths.ox.ac.uk}

\begin{abstract}

This paper is an introduction to the use of perverse sheaves with
positive characteristic coefficients in modular representation theory.
In the first part, we survey results relating singularities in finite
and affine Schubert varieties and nilpotent cones to modular
representations of reductive groups and their Weyl groups.
The second part is a brief introduction to the theory of perverse
sheaves with an emphasis on the case of positive characteristic and
integral coefficients.
In the final part, we provide some explicit examples of stalks
of intersection cohomology complexes with integral or
positive characteristic coefficients in nilpotent cones, mostly in
type $A$. Some of these computations might be new.

\end{abstract}

\maketitle

\section*{Introduction}
\label{sec:introduction}

Representation theory has a very different flavour in positive
characteristic. When one studies the category of representations of a
finite group or a reductive group over a field of characteristic 0
(e.g. $\CM$), one of the first observations to be made is that this
category is semi-simple, meaning that every representation is
isomorphic to a direct sum of irreducible representations.  This
fundamental fact helps answer many basic questions, e.g. the
dimensions of simple modules, character formulae, and tensor product
multiplicities.  However, when one considers representations over
fields of positive characteristic (often referred to as ``modular''
representations) the resulting categories are generally not
semi-simple.  This makes their study considerably more complicated and
in many cases even basic questions remain unanswered.\footnote{
For an introduction to the modular representation theory of finite
groups we recommend the third part of~\cite{Serre}, and for that of
reductive groups,~\cite{Jan}.  }

It turns out that some questions in representation theory
have geometric counterparts. The connection is obtained via the
category of perverse sheaves, a certain category that may be
associated to an algebraic variety and whose structure reflects the
geometry of the underlying variety and its subvarieties. The category
of perverse sheaves depends on a choice of coefficient field and, as
in representation theory, different choices of coefficient field can
yield very different categories.

Since the introduction of perverse sheaves it has been realised that many
phenomena in Lie theory can be explained in terms of categories of perverse sheaves
and their simple objects --- intersection cohomology complexes. In studying
representations of reductive groups and related objects, singular
varieties arise naturally (Schubert varieties and their
generalizations, nilpotent varieties, quiver varieties\dots).  It
turns out that the invariants of these singularities often carry
representation theoretic information. For an impressive list of such
applications, see \cite{LuICM}. This includes constructing
representations, computing their characters, and constructing nice
bases for them.

However, most of these applications use a field $k$ of characteristic
zero for coefficients. In this paper, we want to give the reader a
flavour for perverse sheaves and intersection cohomology with positive characteristic
coefficients.

In the first section of this article we survey three connections between 
modular representation theory and perverse sheaves.

The geometric setting for the first result --- known as the geometric
Satake theorem --- is a space (in fact an ``ind-scheme'') associated to
a complex reductive group $G$. This space, $G(\CM((t)))/G(\CM[[t]])$,
commonly referred to as the affine Grassmannian is a homogeneous space
for the algebraic loop group $G(\CM((t)))$.  Under the action of
$G(\CM[[t]])$, it breaks up as a union of infinitely many
finite-dimensional orbits. Theorems of Lusztig~\cite{Lu2},
Ginzburg~\cite{Gi}, Beilinson-Drinfeld~\cite{BD}, and
Mirkovi\'c-Vilonen~\cite{MV} explain that encoded in the geometry of
the affine Grassmannian and its orbit closures is the algebraic
representation theory over any field (and even over the integers) of
the split form of the reductive group $G^\vee$ with root data dual to
that of $G$, also known as the Langlands dual group.

The second family of results that we discuss involves the geometry of
the finite flag variety $G/B$ where again $G$ is a complex reductive
group, and a generalization of it closely related to the affine
Grassmannian known as the affine flag variety $G(\CM((t)))/\IC$. We
describe theorems of Soergel~\cite{Soe} and Fiebig~\cite{Fie1, Fie2,
Fie3, Fie4} which show that the geometry of these spaces can be used to
understand the modular representation theory of the
Langlands dual group $G^\vee_k$ for $k$ a field of characteristic larger
than the Coxeter number of $G^{\vee}_k$.  
In doing so, Fiebig is able to give a new proof of
the celebrated Lusztig conjecture with an explicit bound on the
characteristic.

The third theorem to be discussed is centered around the geometry of
the variety $\NC$ of nilpotent elements of a Lie algebra $\gg$, known
as the nilpotent cone. The nilpotent cone has a natural resolution and,
in 1976, Springer~\cite{SPRTRIG} showed that the Weyl group acts on the $\ell$-adic
cohomology of the fibers of this resolution.  He showed moreover that
from this collection of representations one could recover all of the
irreducible $\ell$-adic representations and that they came with a
natural labelling by a nilpotent adjoint orbit with an irreducible
$G$-equivariant local system. This groundbreaking discovery was
followed by a series of related constructions, one of which, based on the
Fourier-Deligne transform, has recently been used by the first author~\cite{Ju-thesis}
to establish a modular version of the Springer correspondence.

 The second goal of this
article which occupies the second and third sections is to provide 
an introduction to ``modular'' perverse sheaves,
in other words perverse sheaves with coefficients in a field of
positive characteristic. We begin by recalling the theory of perverse 
sheaves, highlighting the
differences between characteristic zero and characteristic $p$, and
also the case of integer coefficients. We treat in detail the case of
the nilpotent cone of ${\sl}_2$.

In the last part, we treat more examples in nilpotent cones. We
calculate all the IC stalks in all characteristics $\neq 3$ for the nilpotent
cone of ${\sl}_3$, and all the IC stalks in all characteristics $\neq 2$ for
the subvariety of the nilpotent cone of ${\sl}_4$ consisting of the
matrices which square to zero. Before that, we recall how to deal with
simple and minimal singularities in type $A$, for two reasons:
we need them for the three-strata calculations, and they can be dealt
with more easily than for arbitrary type (which was done in
\cite{cohmin,decperv}). As a complement, we give a similar direct approach for
a minimal singularity in the nilpotent cone of $\sp_{2n}$.

The two first parts partly correspond to the talks given by the first
and third author during the summer school, whose titles were
``Intersection cohomology in positive characteristic I, II''. 
The third part contains calculations that the three of us did together
while in Grenoble. These calculations were the first non-trivial
examples involving three strata that we were able to do.

It is a pleasure to thank Michel Brion for organizing this conference
and allowing two of us to speak, and the three of us to meet.  We
would like to thank him, as well as Alberto Arabia, Peter Fiebig, Joel Kamnitzer and
Wolfgang Soergel for very valuable discussions.  The second author
would also like to acknowledge the mathematics department at the
University of Texas at Austin and his advisor, David Ben-Zvi, for
partially funding the travel allowing him to attend this conference
and meet his fellow coauthors.

\setcounter{tocdepth}{2}
\tableofcontents

\section{Motivation}

Perverse sheaves with coefficients in positive characteristic appear
in a number of different contexts as geometrically encoding certain
parts of modular representation theory.  This section will provide a
survey of three examples of this phenomenon: the geometric Satake
theorem, the work of Soergel and Fiebig on Lusztig's conjecture,
and the modular Springer correspondence. The corresponding geometry 
for the three pictures will be respectively
affine Grassmannians, finite and affine flag varieties, and nilpotent
cones.

Throughout, $G_\ZM$ will denote a split connected reductive group
scheme over $\ZM$. Given a commutative ring $k$, we denote by
$G_k$ the split reductive group scheme over $k$ obtained by extension
of scalars
\[
G_k = \Spec k \times_{\Spec \ZM} G_\ZM.
\] 
In Sections \ref{subsec:geomsatake} and
\ref{subsec:SoergelFiebig} we will consider 
$G_\CM$, while in Section \ref{subsec:modSpringer}, we will consider 
$G_{\FM_q}$.

We fix a split maximal torus $T_\ZM$ in $G_\ZM$. We denote by
\[
(X^*(T_\ZM),R^*,X_*(T_\ZM),R_*)
\]
the root datum of $(G_\ZM,T_\ZM)$.
We denote by $(G_\ZM^\vee, T_\ZM^\vee)$ the pair associated to the
dual root datum. Thus $G_\ZM^\vee$ is the Langlands dual group. 
In Subsections \ref{subsec:geomsatake} and
\ref{subsec:SoergelFiebig}, we will consider representations of
$G^\vee_k = \Spec k \times_{\Spec \ZM} G^\vee_\ZM$,
where $k$ can be, for example, a field of characteristic $p$.
We have $X^*(T_\ZM^\vee) = X_*(T_\ZM)$ and $X_*(T_\ZM^\vee) = X^*(T_\ZM^\vee)$.

We also fix a Borel subgroup $B_\ZM$ of $G_\ZM$ containing
$T_\ZM$. This determines a Borel subgroup $B_\ZM^\vee$ of $G_\ZM^\vee$
containing $T_\ZM^\vee$. This also determines bases of simple roots
$\D \subset R^*$ and $\D^\vee \subset R_*$. It will be convenient to
choose $\D_* := -\D^\vee$ as a basis for $R_*$ instead, so that the
coroots corresponding to $B_\ZM^\vee$ are the negative coroots
$R_*^- = - R_*^+$.

\subsection{The geometric Satake theorem}
\label{subsec:geomsatake}

In this subsection and the next one, to simplify the notation, we will
identify the group schemes
$G_\CM \supset B_\CM \supset T_\CM$ with their groups of $\CM$-points
$G \supset B \supset T$.

We denote by $\KC = \CM((t))$ the field of Laurent series 
and by $\OC = \CM[[t]]$ the ring of
Taylor series.  The affine (or loop) Grassmannian $\Gr = \Gr^G$ is the
homogeneous space $G(\KC)/G(\OC)$. It has the structure of an
ind-scheme. In what follows we will attempt to sketch a rough outline
of this space and then briefly explain how perverse sheaves on it are
related to the representation theory of $G^\vee_k$, where $k$ is any
commutative ring of finite global dimension. We refer the reader 
to \cite{BD,BL,LS,MV} for more details
and proofs.

We have a natural embedding of
the coweight lattice $X_*(T) = \Hom(\Gm,T)$ into the affine
Grassmannian: each $\l\in X_*(T)$ defines a point $t^\lambda$ of $G(\KC)$
via
\[
\Spec \KC = \Spec \CM((t)) \elem{c} \Gm = \Spec \CM[t,t^{-1}]
\elem{\l} T \elem{i} G
\]
where $c$ comes from the inclusion
$\CM[t,t^{-1}] \hookrightarrow \CM((t))$
and $i : T \to G$ is the natural inclusion,
and hence a point $[t^\lambda]$ in $\Gr = G(\KC)/G(\OC)$.

For example, when $G=GL_n$ and $T$ is the subgroup of diagonal
matrices the elements of $X_*(T)$ consist of $n$-tuples of integers
$\lambda=(\lambda_1, \ldots, \lambda_n)$ and they sit inside of $\Gr$
as the points
\[
\left( \begin{array}{cccc} t^{\lambda_1} & & & \\ &
  t^{\lambda_2} & & \\ & & \ddots & \\ & & & t^{\lambda_n}
 \end{array} \right)\cdot G(\OC)
\]

As in the finite case, one has a
Cartan decomposition of $G(\KC)$
\[
G(\KC)=\bigsqcup_{\lambda\in X_*(T)^+} G(\OC)t^\lambda G(\OC)
\]
where $X_*(T)^+=\{\lambda\in
X_*(T) \mid \forall\a\in\Delta,\ \langle\alpha,\lambda\rangle \geq 0\}$
is the cone of dominant coweights. Thus
the affine Grassmannian is the union of the $G(\OC)$-orbits of the
points $[t^\lambda]$ for $\lambda\in X_*(T)^+$.

Another important feature of the affine Grassmannian is a special
$\CM^*$-action. As in the topological loop group, there is a notion of
``loop rotation''. In our case, this rotation comes in the form of
$\CM^*$ acting on $G(\KC)$ by composing a $\KC$-point of $G$ with the
automorphism of $\Spec \KC$ induced by the action of $\CM^*$ on
itself.  More naively, this means replacing $t$ by $zt$.  This action
clearly preserves the subgroup $G(\OC)$ and thus gives a well-defined
action on the quotient $\Gr$.

It is useful to get a sense of the geometry of the $G(\OC)$-orbits.  As we saw
above, these orbits are labelled by the dominant coweights. We begin
by studying a subvariety of the $G(\OC)$-orbits. For
$\lambda\in X_*(T)^+$, consider the $G$-orbit of $[t^\lambda]$.  It
turns out that for such a dominant coweight, the point $[t^\lambda]$ is
fixed by a Borel subgroup. Thus, the $G$-orbit is a (partial) flag
variety. In fact, the stabilizer of $[t^\lambda]$ in $G$ is a
parabolic subgroup $P_\lambda$ with Levi factor corresponding to the
roots $\alpha \in \Delta$ such that
$\langle\alpha,\lambda\rangle=0$. We conclude that $G\cdot
[t^\lambda]$ is isomorphic to the (partial) flag variety
$G/P_\lambda$. It is an easy exercise to check these claims for
$G=GL_n$.

The points $[t^\lambda]$ of $\Gr$ are in fact fixed by loop rotation.  To 
see this, note that the reduced affine Grassmannian for $T$, $T(\KC)/T(\OC) 
\cong X_*(T)$ is discrete, embeds in $\Gr^G$ as the set of points 
$[t^\lambda]$ for $\lambda\in X_*(T)$, and is preserved by loop rotation. As 
$\CM^*$ is connected and the subset $[t^\lambda]$ discrete, each such 
point is fixed under loop rotation.  But of course the group
$G\subset G(\OC)$ is certainly fixed by loop rotation, thus the
$G$-orbit $G\cdot [t^\lambda]$ is fixed under loop rotation as well. 

Not only are these $G$-orbits $G\cdot [t^\lambda]$ fixed,
they form precisely the fixed point set of the action of loop rotation
on the affine Grassmannian. Moreover, the $G(\OC)$-orbit 
$G(\OC)\cdot [t^\lambda]$ is a vector bundle over $G\cdot
[t^\lambda]\cong G/P_\lambda$.
A proof of this statement involves considering the highest congruence
subgroup of $G(\OC)$, defined as the preimage $\mathrm{ev}_0^{-1}(1)$ of $1$ under the 
evaluation map $\mathrm{ev}_0 : G(\OC)\to G, t \mapsto 0$. One can check that the
orbit $\mathrm{ev}_0^{-1}(1) \cdot g \cdot [t^\lambda]$ is an affine
space for any $g\in G$ and is isomorphic to a vector space on which
loop rotation acts linearly by contracting characters.  Combining this
with the fact that $G(\OC)=\mathrm{ev}^{-1}_0(1)\cdot G$, 
the claim follows. A corollary of this remark is that the $G(\OC)$-orbits are simply-connected.

For $\lambda$ and $\mu$ dominant, the orbit $G(\OC)\cdot 
[t^\l]$ is of dimension $2\rho(\l)$ (here $\rho$ is half the sum 
of the positive roots) and is contained in $\overline{G(\OC)\cdot 
[t^\mu]}$ if and only if $\lambda-\mu$ is a sum of positive coroots. 

As a concrete example, it is instructive to consider the case
$G=PSL_2$. Choose $T$ to be the subgroup of diagonal matrices (up to
scale) and $B$ the upper triangular matrices (again, up to scale). The
torus $T$ is one dimensional, so the lattice $X_*(T)$ is isomorphic to
a copy of the integers and $X_*(T)^+$ to the non-negative ones. Thus
the $G(\OC)$-orbits are labelled by the non-negative integers. The
parabolic subgroup corresponding to any positive number is the Borel
subgroup $B$, and that correpsonding to the trivial weight is the
whole group $PSL_2$. Thus the affine Grassmannian for $PSL_2$ is a
union of a point and a collection of vector bundles over
$\PM^1$. Considering the remark of the previous paragraph, as the
coroot lattice for $PSL_2$ is a subgroup of index two in the coweight
lattice, the affine Grassmannian consists of two connected components.

\begin{remark} One way to see some of the geometry is through the moment 
map with respect to the action of the torus $T$ extended by loop 
rotation.  This idea, from the differential point of view and in 
slightly different language, can be found in ~\cite{ap}.  Yet another 
picture of the affine Grassmannian is provided by the spherical 
building, whose vertices are the $\FM_q$-points of the affine 
Grassmannian where the Laurent and Taylor series are defined over 
$\overline\FM_p$ instead of $\CM$. This picture in the rank one case can 
be found in chapter 2 of~\cite{s}, although the affine Grassmannian is 
not mentioned explicitly. \end{remark}

Beginning with pioneering work of Lusztig, it was understood that the
geometry of the affine Grassmannian is closely related to the
representation theory of the Langlands dual group $G^\vee$, i.e. the reductive group with dual root data.  In particular,
Lusztig showed~\cite{Lu2} that the local intersection cohomology (with
complex coefficients) of the $G(\OC)$-orbits was a
refinement of the weight multiplicities of the corresponding
representation of $G^\vee$.

This connection was further developed by Ginzburg~\cite{Gi} (see also
Beilinson-Drinfeld~\cite{BD}), who noted that the category of
$G(\OC)$-equivariant perverse sheaves (with $\CM$-coefficients) carried a
convolution product and using Tannakian formalism on the total
cohomology functor that the category was tensor equivalent to the
category of representations of the Langlands dual group $G^\vee$. In
other words:

\begin{thm}
There is an equivalence of tensor categories:
\[
(\Perv_{G(\OC)}(\Gr^G;\CM),*) \simeq 
(\Rep(G^\vee_\CM),\otimes).
\]
\end{thm}
  
This result can be interpreted as a categorification of the much
earlier work of Satake~\cite{Sat} which identified the algebra of spherical
functions (approximately bi-$G(\OC)$-invariant functions of $G(\KC)$)
with $W$-invariant functions on the coweight lattice,
$\CM[\Lambda]^W$.

\begin{remark} In the case of the affine Grassmannian, the category of
  $G(\OC)$-equivariant perverse sheaves is equivalent to the category
  of perverse sheaves constructible with respect to the $G(\OC)$-orbit
  stratification. For a proof, see the appendix to~\cite{MV}.
\end{remark}

It was understood by Beilinson and Drinfeld~\cite{BD} that the affine
Grassmannian described above should be thought of as associated to a
point on an algebraic curve.  Understood as such, there is a natural
global analogue of the affine Grassmannian living over a configuration
space of points on a curve. Using this Beilinson-Drinfeld
Grassmannian, one can produce a natural commutativity constraint for
the convolution product by identifying it as a ``fusion product''.
>From this point of view, the geometric Satake theorem becomes
identified with the local geometric Langlands conjecture.

For the remainder of this paragraph let $k$ be a Noetherian
commutative ring of finite global
dimension. Mirkovi{\'c}-Vilonen~\cite{MV} generalized and rigidified
the picture further by producing the analogue of the weight functors
for perverse sheaves with coefficients in an arbitrary $k$.  Consider
the functor $F_\nu:\Perv_{G(\OC)}(\Gr;k)\to k\modules$ for each $\nu\in
X^*$ which takes compactly supported cohomology along the
$N(\KC)$-orbit containing $[t^\nu]$. They prove that these cohomology
groups vanish outside of degree $2\rho(\nu)$ and that the functors are
exact.  Summing over all $\nu$, they prove that there is a natural
equivalence of functors
\[
H^* \cong \bigoplus_{\nu\in X_*(T)} F_\nu : \Perv_{G(\OC)}(\Gr;k)\to k\modules.
\]
This more refined fiber functor together with some delicate arithmetic
work allowed them to prove that geometric Satake is true for any such
$k$, meaning the category of $G(\OC)$-equivariant perverse sheaves
with $k$-coefficients is tensor equivalent to the category of
representations of the split form of the Langland dual group $G^\vee$
over $k$. In other words,

\begin{thm}
There is an equivalence of tensor categories:
\[
(\Perv_{G(\OC)}(\Gr^G;k),*) \simeq (\Rep(G^\vee_k),\otimes).
\]
\end{thm}

\subsection{Finite and affine flag varieties}
\label{subsec:SoergelFiebig}

In this subsection we give an overview of work of Soergel \cite{Soe}
and Fiebig \cite{Fie1, Fie2, Fie3} relating the rational
representation theory of reductive algebraic groups over a field $k$
of positive characteristic to sheaves on complex Schubert varieties
with coefficients in $k$.

Fix a field $k$ of characteristic $p$. Recall that $G$ denotes a 
reductive algebraic group 
over $\CM$ and that $G^{\vee}_k$ is the split reductive algebraic 
group over a field $k$ with root datum dual to that of $G$. In this
section we assume that $G$ is connected, simple and adjoint. It follows that 
$G^{\vee}_k$ is simply connected.

The previous section explained how one may give a geometric
construction of the entire category of representations of $G^{\vee}_k$
in terms of $G$. The constructions which follow establish a relation
between blocks (certain subcategories) of representations of
$G^{\vee}_k$ and sheaves on (affine) Schubert varieties associated to
$G$.

In order to explain this we need to recall some standard facts from
representation theory which one may find in \cite{Jan}. Recall that 
we have
also fixed a Borel subgroup and maximal torus
$G^{\vee}_k \supset B^{\vee}_k \supset T^{\vee}_k$,
that we write $R_*$ and $R_*^+$ for the roots and positive roots of
$(G^\vee_k,T^\vee_k)$ respectively,
chosen so that $-R_*^+$ are the roots 
determined by $B^\vee_k$.  By duality we may identify $X_*(T)$ and 
$X^*(T^{\vee}_k)$.
We denote by $\Rep G^{\vee}_k$ the category of all finite
dimensional rational representations of $G^{\vee}_k$.

To each weight $\lambda \in X^*(T^{\vee}_k)$ one may associate a $G^{\vee}_k$-equivariant
line bundle $\OC(\lambda)$ on $G^{\vee}_k/B^{\vee}_k$. Its global sections
\begin{equation*}
H^0(\lambda) = H^0(G^{\vee}_k/B^{\vee}_k, \OC(\lambda))
\end{equation*}
contain a unique simple subrepresentation $L(\lambda)$, and all simple
$G^{\vee}_k$-modules arise in this way. The module $H^0(\lambda)$ is
non-zero if and only if $\lambda$ is dominant.

It is known that the characters of $H^0(\lambda)$ are given by the
Weyl character formula and that $L(\mu)$ can only occur as a
composition factor in $H^0(\lambda)$ if $\lambda - \mu \in \NM R_*^+$. It
follows that in order to determine the characters of the simple
$G^{\vee}_k$-modules, it is enough to determine, for all dominant 
$\lambda, \mu \in X^*(T^\vee_k)$, the multiplicities:
\begin{equation}
\label{h0-mult}
[H^0(\lambda): L(\mu)] \in \NM.
\end{equation}

In fact, many of these multiplicities are zero. Recall that the Weyl
group $W$ acts on $X^*(T^{\vee}_k)$ and we may consider the ``dot action'' given
by
\begin{equation*}
w \cdot \lambda = w(\lambda + \rho) - \rho
\end{equation*}
where $\rho$ denotes the half-sum of the positive roots. Denote by
$\widehat{W}$ the subgroup of all affine transformations of $X^*(T^\vee_k)$
generated by $( w \cdot )$ for $w \in W$ and $( + \mu)$ for all $\mu
\in p \ZM R_*$. As an abstract group this is isomorphic to the affine Weyl
group associated to the root system of $G^{\vee}_k$. The category of rational
representations of $G^{\vee}_k$ decomposes into blocks\footnote{A
  family of full subcategories $\CC_{\Omega}$ of a category $\CC$
  yields a block decomposition (written $\CC = \bigoplus
  \CC_{\Omega}$) if every object $M$ in $\CC$ is isomorphic to a direct
  sum of objects $M_\Omega\in\CC_{\Omega}$ and there are no morphisms
  between objects of $\CC_{\Omega}$ and objects of
  $\CC_{\Omega^{\prime}}$ if $\Omega \ne \Omega^{\prime}$.
}
\begin{equation*}
\Rep G^{\vee}_k = \bigoplus_{\Omega} \Rep_{\Omega} G^{\vee}_k
\end{equation*}
where $\Omega$ runs over the orbits of $\widehat{W}$ on
$X^*(T^{\vee}_k)$ and $\Rep_{\Omega} G^{\vee}_k$ denotes the full
subcategory of $\Rep G^{\vee}_k$ whose objects are those
representations all of whose simple factors are isomorphic to
$L(\lambda)$ for $\lambda \in \Omega$.

Assume from now on that $p > h$, where $h$ denotes the Coxeter number
of the root system of $G^{\vee}_k$. Then the ``translation principle'' allows
one to conclude that it is enough to know the characters of the simple
modules in $\Rep_{\Omega} G^{\vee}_k$, where $\Omega = \widehat{W} \cdot 0$. This
fact, combined with the Steinberg tensor product theorem, allows one
to reduce the problem to calculating the multiplicities
\begin{equation} \label{h0dot-mult}
[H^0(x \cdot 0): L(y \cdot 0)] \in \NM
\end{equation}
where $x \cdot 0$ and $y \cdot 0$ lie in the ``fundamental box'':
\begin{equation*}
I = \{ \lambda \in X^*(T^{\vee}_k) \; | \; \langle \alpha^{\vee}, \lambda \rangle
< p \text{ for all } \alpha^{\vee} \text{ simple} \}.
\end{equation*}
A celebrated conjecture 
of Lusztig 
expresses the multiplicity in
(\ref{h0dot-mult}) for $x \cdot 0$ and $y \cdot 0$ lying in $I$ in
terms of certain Kazhdan-Lusztig polynomials evaluated at 1.
Kazhdan-Lusztig polynomials are the Poincaré polynomials of the local
intersection cohomology of Schubert varieties (here in the affine
case) \cite{KL2}, and they can be defined through an inductive
 combinatorial procedure (which provided their original definition in
\cite{KL1}). This conjecture is known to hold for almost all $p$ 
by work of Andersen, Jantzen and Soergel \cite{AJS}.

We now return to geometry. Recall that $G$ is a reductive algebraic 
group whose root system
is dual to that of $G^{\vee}_k$. We identify the Weyl groups of $G$
and $G^{\vee}_k$ in the obvious way. For each simple reflection $s \in
S$ we may associate a minimal parabolic $B \subset P_s \subset G$ and
we have a projection map
\begin{equation*}
\pi_s : G/B \to G/P_s.
\end{equation*}
Let $D^b_c(G/B, k)$ denote the bounded derived category of constructible sheaves of
$k$-vector spaces on $G/B$. In \cite{Soe}, Soergel considers the
category $\KC$ defined to be the smallest additive subcategory of
$D^b(G/B; k)$ such that:
\begin{enumerate}
\item the skyscraper sheaf on $B/B \in G/B$ is in $\KC$;
\item if $\FS \in \KC$ then so is $\pi_s^* {\pi_s}_* \FS$;
\item if $\FS\in \KC$ then so is any object isomorphic to a shift of a
  direct summand of $\FS$.
\end{enumerate}
If $k$ were of characteristic 0 then one could use the decomposition
theorem to show that any indecomposable object in $\KC$ isomorphic to
a (shift of) an intersection cohomology complex of the closure of a
$B$-orbit on $G/B$. However, as $k$ is of positive
characteristic this is not necessarily the case. Somewhat
surprisingly, for each $x \in W$, it is still true that there exists
up to isomorphism a unique indecomposable object $\FS_x \in \KC$ 
supported on $\ov{ BxB/B}$ and such
that $(\FS_x)_{B x B / B} \simeq \uk_{BxB /B}[\ell(w)]$. Each $\FS_x$ is
self-dual and any indecomposable object in $\KC$ is
isomorphic to $\FS_x[m]$ for some $x \in W$ and $m \in \ZM$.

Soergel goes on to establish a connection between $\KC$ and the
representation theory of $G^{\vee}_k$ as follows. 
Let
\[ 
\rho = \frac{1}{2} \sum_{\alpha \in R_*^+} \alpha \in X^*(T^\vee_k)
\]
and $st = (p - 1) \rho$ be
the Steinberg weight. He shows:

\begin{thm}[\cite{Soe}, Theorem 1.2] 
With $k$ as above, for $x, y \in W$, one has
\begin{equation*}
[H^0( st + x \rho): L(st + y \rho)] = \dim (\FS_y)_x
\end{equation*}
where $\dim(\FS_y)_x$ denotes the total dimension of the cohomology of
the stalk
$(\FS_y)_x$.
\end{thm}

A disadvantage of the above approach is that it offers a geometric
interpretation for only a small part of the representation theory of
$G^{\vee}_k$. In recent work Fiebig has developed a more complete (and
necessarily more complicated) picture.

Before we describe Fiebig's work we recall a construction of
$T$-equivariant intersection cohomology due to Braden and 
MacPherson in \cite{BrM}. Let $T \simeq
(\CM^*)^n$ be an algebraic torus and $X$ a complex $T$-variety with
finitely many zero- and one-dimensional orbits. To this situation one 
may associate a labelled graph called the
``moment graph'', which encodes the structure of the zero- and
one-dimensional orbits (see the notes of Jantzen from this
conference).
Under some additional assumptions on $X$ (the most important being a
$T$-invariant stratification into affine spaces) Braden and MacPherson
describe a method to calculate the $T$-equivariant intersection
cohomology of $X$ with coefficients in $\QM$. This involves the 
inductive construction of
 a ``sheaf'' $M(X, \QM)$ on the moment graph of $X$; the
equivariant intersection cohomology is then obtained by taking
``global sections'' of this sheaf.

As Fiebig points out, Braden and MacPherson's construction makes sense
over any field $k$ and produces a sheaf on the moment graph $M(X, k)$;
however it is not clear if what one obtains in this way has anything
to do with the intersection cohomology with coefficients in $k$.

Let $G$ and $B$ be as above and consider the affine flag variety $
G((t)) / \IC$ where $\IC$ denotes the Iwahori subgroup, defined as the
preimage of $B$ under the evaluation map $G(\OC)\to G,t\mapsto 0$.
As with the loop Grassmannian, $
G((t)) / \IC$ may be given the structure of an ind-scheme.  The
$\IC$-orbits on $G((t))/\IC$ are affine spaces parametrized by
$\widehat{W}$. Recall that earlier we defined the fundamental box $I
\subset X^*(T^\vee_k)$. Define
\begin{equation*}
\widehat{W}^\res = \{ w \in \widehat{W} \; | \; - w \cdot 0 \in I \}.
\end{equation*}
Fiebig shows:

\begin{thm}
If the stalks of the sheaf on the moment graph $M(\overline{\IC x \IC/
  \IC}, k)$ are given by Kazhdan-Lusztig polynomials for all $x \in
\widehat{W}^\res$ then Lusztig's conjecture holds for
representations of $G^{\vee}_k$.
\end{thm}

For a more precise statement we refer the reader to \cite{Fie3}. This
result enables Fiebig to give a new proof of the Lusztig's conjecture
in almost all characteristics. Actually, Fiebig is able to give an
explicit bound, although it is still very big \cite{Fie4}. Moreover
he is able to prove the multiplicity one case of the conjecture for all
characteristics greater than the Coxeter number \cite{Fie1}.

It is expected that one may obtain the sheaf on the moment graph
$M(\overline{\IC x \IC/ \IC}, k)$ by applying a functor similar to that
considered by Braden and MacPherson to a sheaf $\FS_x \in
D^b_{G[[t]]}(G((t))/\IC; k)$. (This sheaf should be analogous to the
indecomposable sheaves considered by Soergel.) If this is the case
then Fiebig's theorem asserts that Lusztig's conjecture would follow
from a certain version of the decomposition theorem with coefficients in
$k$.

\subsection{The modular Springer correspondence}
\label{subsec:modSpringer}

In this subsection, the base field is not $\CM$ but $\FM_q$, where $q$
is a power of some prime $p$. Perverse sheaves still make sense in
this context, using the étale topology \cite{BBD}.
Now $G$ will be $G_{\FM_q}$, which we identify its set of
$\ov\FM_p$-points, endowed with a Frobenius endomorphism $F$.
We denote by $\gg$ its Lie algebra, and by $W$
its Weyl group. For simplicity, we assume that $p$ is very good for
$G$, so that the Killing form provides a non-degenerate $G$-invariant
symmetric bilinear form on $\gg$. Thus we can identify $\gg$ with its
dual $\gg'$.

In 1976, Springer established a link between the \emph{ordinary}
(that is, characteristic zero) representations of $W$, and the
nilpotent cone $\NC$ of $\gg$ \cite{SPRTRIG}. More precisely, he
constructed the irreducible representations of $W$ in the top
cohomology (with $\ell$-adic coefficients, $\ell$ being a prime
different from $p$), of some varieties attached to the different
nilpotent orbits, the Springer fibers. To each irreducible
representation of $W$ corresponds a nilpotent orbit and an irreducible
$G$-equivariant local system on this orbit.

The modular Springer correspondence \cite{Ju-thesis} establishes such
a link for \emph{modular} representations of $W$, over a field of
characteristic $\ell$. The modular irreducible representations of $W$
are still largely unknown, for example if $W$ is a symmetric group of
large rank. One would like to know their characters, and this is
equivalent to determining the entries in the so-called decomposition
matrix, relating ordinary and modular irreducible characters. Using
the modular Springer correspondence, one can show that this
decomposition matrix can be seen as a submatrix of a decomposition
matrix for $G$-equivariant perverse sheaves on the nilpotent cone.
As a result, just as the geometric Satake isomorphism implies that the
modular representation theory of reductive groups is encoded in the
singularities of the complex affine Grassmannian of the dual group, one can
say that the modular representation theory of the Weyl group of a Lie
algebra is encoded in the singularities of its nilpotent cone.

We fix a Borel subgroup $B$ of $G$, with Lie algebra $\bg$. We denote
by $U$ the unipotent radical of $B$, and by $\ug$ the Lie algebra of
$U$. Then $\ug$ is the orthogonal of $\bg$.
The group $G$ acts transitively on the set of Borel subalgebras of
$\gg$, and the normalizer of the Borel subalgebra $\bg$ is $B$, so the
flag variety $\BC := G/B$ can be identified with the set of all Borel
subalgebras of $\gg$. It is a smooth projective variety.

One then defines $\tilde\NC :=\{(x,gB)\in \NC\times G/B \mid x\in
\Ad(g)\bg\}$. One can check that the second projection makes
it a $G$-equivariant vector bundle over $\BC = G/B$, and that
we have $G$-equivariant isomorphisms:
\[
\tilde \NC \simeq G\times^B \ug
\simeq G\times^B \bg^\perp
\simeq G\times^B (\gg/\bg)^*
\simeq T^*(G/B) = T^*\BC
\]
where the first isomorphism is given by $(y,gB) \mapsto g*\Ad(g^{-1})y$.
The first projection gives a resolution $\pi_\NC : \tilde\NC \to \NC$
of the nilpotent cone, called the Springer resolution.

Springer constructed an action of the Weyl group $W$ on the
$\ell$-adic cohomology
of the fibers $\BC_x := \pi_\NC^{-1}(x)$ of this resolution, which are
called Springer fibers. These are connected projective varieties, which
are usually singular. All their irreducible components have the same dimension $d_x :=
\frac 1 2 \codim_\NC G.x$. In particular, to each adjoint orbit of
$\NC$ one can associate the representation of $W$ on the top cohomology of
the corresponding Springer fiber. In type A, this is in fact a
bijection between the nilpotent orbits and the irreducible
representations of $\SG_n$. Note that both are parametrized by the set
of all partitions of the integer $n$. It turns out that the bijection
is given by the conjugation of partitions.

More generally, for a point $x\in \gg$, let $C_G(x)$ denote the
centralizer of $x$ in $G$ and $A_G(x)=C_G(x)/C_G(x)^o$ its component
group; for $G = GL_n$ these groups are trivial. As $\pi_\NC$ is
$G$-equivariant, the centralizer $C_G(x)$ acts on the fiber $\BC_x$,
and the group $A_G(x)$ acts on the cohomology of $\BC_x$. This action
commutes with the action of $W$. Note that the action of $A_G(x)$ on
the top cohomology of $\BC_x$ is just the permutation representation
of $A_G(x)$ on the set of irreducible components of $\BC_x$. It turns out that
$H^{2d_x}(\BC_x)$ is an irreducible $W \times A_G(x)$-module.
We can decompose it into $A_G(x)$-isotypic components:
\[
H^{2d_x}(\BC_x) = \bigoplus_\rho \rho \otimes V_{x,\rho}
\]
where $\rho$ runs over all irreducible representations of $A_G(x)$ such
that the $\rho$-isotypic component of $H^{2d_x}(\BC_x)$ is non-zero,
and $V_{x,\rho}$ is a well-defined irreducible representation of
$W$. Springer showed that the $V_{x,\rho}$, for $x$ running over a set
of representatives of the nilpotent orbits, form a complete collection
of irreducible representations of the Weyl group $W$. So to each
irreducible representation of $W$ we can assign a pair
$(x,\rho)$. This is the Springer correspondence.

Later a number of related constructions were obtained by other
mathematicians, as in \cite{SLO1,KL}. It turns out that the pairs
$(x,\rho)$ also parametrize the simple $G$-equivariant perverse
sheaves on the nilpotent cone. Let us denote them by
$\ic(x,\rho)$. Lusztig and Borho-McPherson gave a construction using
perverse sheaves \cite{lu,BM2}. Note, however, that all these
approaches give a parametrization which differs from the original one
obtained by Springer by tensoring with the sign character of $W$. On
the other hand, other approaches, using some kind of Fourier
transform, give the same parametrization as the one by Springer. One
can use a Fourier transform for $\DC$-modules \cite{HK} if the base
field is $\CM$, or a Fourier-Deligne transform in our context where
the base field is $\FM_q$ \cite{BRY}, using perverse $\KM$-sheaves,
where $\KM$ is a finite extension of $\QM_\ell$. One advantage of the
latter approach is that it still makes sense over finite extensions
$\OM$ (resp. $\FM$), of $\ZM_\ell$ (resp. $\FM_\ell$), so that it is
possible to define a modular Springer correspondence \cite{Ju-thesis}.

The Fourier-Deligne transform is an equivalence of derived categories
of constructible sheaves between $V$ and $V'$, where $V$ is a vector
bundle $\xi:V \to S$ of constant rank $r$ over a scheme $S$ of finite
type over $k$, and $\xi':V' \to S$ is its dual. In particular, if $S =
\Spec \FM_q$, then $V$ is just an $\FM_q$-vector space. For $\gg \to
\Spec \FM_q$, we get an auto-equivalence of $D^b_c(\gg,\KM)$, since we
identify $\gg$ with its dual. Though we will use the Fourier-Deligne
transform in the particular case $S = \Spec \FM_q$, at some point we
will also need to use the relative version, with $S = \BC$.

This equivalence is a sheaf-theoretic version of the ordinary
Fourier transform for functions on $\RM^n$.  Recall that the ordinary
Fourier transform of a function $f : \RM^n \to \CM$
is the function $\hat f : (\RM^n)^* \to \CM$ given by the formula:
\[
\hat{f} (\zeta) = \int_{\RM^n} f(x)e^{-2\pi i \langle x, \zeta \rangle}dx.
\]
In other words, the Fourier transform takes a function $f$ on $\RM^n$, and
\begin{itemize}
\item pulls it back by the first projection to the product $\RM^n\times(\RM^n)^*$;
\item multiplies it by the pull-back of the exponential function
  $t\mapsto e^{-2\pi i t}$ via the evaluation map $\RM^n \times (\RM^n)^* \to
  \RM$;
\item pushes it forward to $(\RM^n)^*$ by integrating along the fibres
  of the second projection.
\end{itemize}

In order to mimic this procedure sheaf-theoretically, we need to find
a replacement for the exponential function. This role is played by an
Artin-Schreier sheaf, and this is the reason why we use $\FM_q$ as a
base field, instead of a field of characteristic zero.

First, let us define a Fourier transform for $\CM$-valued functions on
$\FM_q^n$. Instead of the exponential character $\RM \to \CM^*$,
$t \mapsto e^{-2\pi i t}$, we have to choose a non-trivial character $\psi$ of
$\FM_q$. For example, we can take the character
$t\mapsto e^{-\frac{2\pi i}p t}$ of $\FM_p$, and compose it with
$\Tr_{\FM_q/\FM_p}$. For $f : \FM_q^n \to \CM$, we set
\[
\hat{f} (\zeta) = 
\sum_{x\in\FM_q^n} f(x)e^{-\frac{2\pi i}p \psi(\langle x, \zeta \rangle)}.
\]

Let $k$ denote either $\KM$, $\OM$ or $\FM$ (see above). We can
replace $\CM$-valued functions by $k$-valued functions, as soon as
$k^\times$ contains the $p$th roots of unity, which we assume from now
on.

Let us now define the Fourier-Deligne transform. We consider the
Artin-Schreier covering of the affine line: $\AM^1 \to \AM^1$,
$t\mapsto t-t^q$. It is a Galois finite étale morphism, with Galois
group $\FM_q$. Thus to the character $\psi$ corresponds a local system
on $\AM^1$, which we will denote by $\LC_\psi$. We can pull it back to
$V \times_S V'$ by the canonical pairing $\mu:V \times_S V' \to
\AM^1$.

As with functions, the Fourier-Deligne transform is a
convolution against a kernel. It is defined as
\[
\FC = (\pr')_!(\pr^*(-) \otimes \mu^* \LC_\psi):
D^b_c(V,k) \longto D^b_c(V',k),
\]
where the notation is fixed by the following diagram:
\begin{equation}
\xymatrix{
&V\times_S V' \ar[ld]_{\pr} \ar[rd]^{\pr'} \ar[rr]^{\mu} &&\AM^1\\
V\ar[rd] && V' \ar[ld]\\
& S &
}
\end{equation}

Most properties of the ordinary Fourier transform still hold for the
Fourier-Deligne transform, if translated appropriately:
it interchanges the extension by zero of the constant sheaf on the
zero section with the
constant sheaf in degree $-\dim V$ (and more generally
interchanges the constant sheaf on a sub-bundle with 
the constant sheaf on its annihilator up to a shift), 
it ``squares to the identity'' up to a sign, which implies that it is an
equivalence of triangulated categories, and it behaves well
under base-change. Moreover, it takes perverse sheaves to perverse
sheaves and, restricted to such, is in fact an equivalence of abelian
categories.

We will now briefly describe the Springer correspondence using this
Fourier-Deligne transform. Let us consider the adjoint quotient $\chi
: \gg \to \gg /\!/G$. The morphism $\chi$ is flat and surjective. By
Chevalley's restriction theorem, we have $\gg/\!/G \simeq \tg/W$. We
denote by $\phi : \tg \to \tg/W$ the quotient morphism, which is
finite and surjective. For $t \in \tg$, we will write $\phi(t) = \ov
t$.  One can see the nilpotent cone $\NC$ as $\chi^{-1}(\ov 0)$. On
the other hand, if $t$ is a regular element in $\tg$, then
$\chi^{-1}(\ov t) \simeq G/T$ is smooth. The fibers $\chi^{-1}(\ov
t)$, for $\ov t$ varying in $\tg/W$, interpolate these two extreme
cases. Grothendieck found a way to obtain a resolution for all the
fibers of $\chi$ simultaneously. We introduce
$\tilde \gg := \{(x,gB)\in \gg\times G/B \mid x\in \Ad(g)\bg\}$.
The second projection makes it a $G$-equivariant vector bundle over
$\BC$, isomorphic to $G\times^B\bg$. The first projection defines a
proper surjective morphism $\pi : \tilde \gg \to \gg$.
Then one can form a commutative diagram:
\[
\xymatrix{
\tilde\gg \ar[r]^\pi\ar[d]_\theta &
\gg \ar[d]^\chi
\\
\tg \ar[r]_\phi &
\tg/W
}
\]
where $\theta$ is the map $G\times^B \bg \to
\bg/[\bg,\bg] \simeq \tg$.
This is a smooth surjective morphism. Then for all $t \in \tg$, the
morphism $\pi_t : \theta^{-1}(t) \to \chi^{-1}(\ov t)$ is a resolution
of singularities.

In the case of ${\gl}_n$, one can see $\chi$ as the map taking
a matrix to its eigenvalues (with multiplicities) up to ordering.
On the other hand, the variety $G/B$ can be identified with the
variety of all complete flags
$F_\bullet = (0 \subset F_1 \subset \dots \subset F_n = \CM^n)$.
Then $\theta$ is identified with the map taking a pair $(x,F_\bullet)$ to
$(a_1,\ldots,a_n)$ where $a_i$ is the eigenvalue of $x$ on
$F_i/F_{i-1}$.

Inside of $\gg$, we have the dense open subvariety $\gg_{rs}$ of
regular semi-simple elements. Consider for a moment the case of
${\gl}_n$, then the semi-simple regular elements are those that
are diagonalizable with pair-wise distinct eigenvalues. For such a
matrix consider the set of Borel subalgebras containing it.  This is
equivalent to the set of full flags preserved by the matrix. But the
matrix decomposes $\CM^n$ into a direct sum of $n$ lines with distinct
eigenvalues.  Thus any subspace preserved by the matrix is a sum of
such lines.  It follows that the collection of such flags is a torsor
for the symmetric group on the eigenspaces. More generally, the
restriction $\pi_\rs : \tilde\gg_\rs \to \gg_\rs$ of $\pi$ over
$\gg_\rs$ is a principal Weyl group bundle: one can show
that it is a Galois finite étale morphism with Galois group $W$.

A local system on $\gg_\rs$ can be identified with a representation of
the fundamental group of $\gg_\rs$. By the above, $W$ is a finite
quotient of this fundamental group. Actually, the local systems which
correspond to a representation factoring through $W$ are those whose
pull-back to $\tilde\gg_\rs$ is trivial. The local system
$\pi_{\rs *} \uk_{\tilde\gg_\rs}$ corresponds to the regular
representation of the group algebra $kW$. Its endomorphism algebra is
the group algebra $kW$.

It is known that
$\pi_* \uk_{\tilde \gg}[\dim \gg] = \ic(\gg,\pi_{\rs *}
\uk_{\tilde\gg_\rs})$
because $\pi$ is a small proper morphism (this notion will be recalled
in Subsection \ref{subsec:perv field}). But, in general, the IC
complex is given by an intermediate extension functor
$j_{\rs !*}$ which is fully faithful, where $j_\rs:\gg_\rs \to \gg$
is the open immersion. Thus the endomorphism algebra of
$\pi_* \uk_{\tilde \gg}[\dim \gg]$ is still $kW$.

We have
\[
\FC(\pi_{\NC*}\uk_{\tilde\NC}[\dim \NC])
\simeq \pi_*\uk_{\tilde\gg}[\dim \gg]
\]
(we ignore Tate twists). To see this one uses the descriptions of
$\tilde\NC$ and $\tilde\gg$ as $G\times^B \ug$ and $G\times^B
\bg$ respectively. These are two orthogonal sub-bundles of the trivial
bundle $\gg\times\BC$ over $\BC$, hence they are exchanged by the
Fourier-Deligne transform with base $\BC$. By base change, it follows
that the perverse sheaves
$\pi_{\NC*} \uk_{\tilde\NC}[\dim \NC]$ and $\pi_*\uk_{\tilde\gg}[\dim \gg]$ are
exchanged by the Fourier-Deligne transform on $\gg$ (here the former
is considered as a perverse sheaf on $\gg$ by extension by zero).
As the Fourier-Deligne transform is an equivalence of categories, we
can conclude that the endomorphism algebra of 
$\pi_{\NC*} \uk_{\tilde\NC}[\dim \NC]$ is again $kW$.

Now, assume $k$ is $\KM$ or $\FM$. Given a simple $kW$-module $E$, one
can consider the corresponding local system $\LC_E$ on $\gg_\rs$.
Then one can show that $\FC(\ic(\gg,\LC_E))$ is a simple
$G$-equivariant perverse sheaf supported on $\NC$, thus we can
associate to $E$ a nilpotent orbit and a representation of the
centralizer component group. In this way, one obtains a Springer
correspondence for arbitrary characteristic: 
a map $\Psi_k$ from $\Irr kW$ to the set $\PG_k$ 
of such pairs consisting of an orbit and representation of
the component group (for characteristic zero coefficients, see
\cite{BRY}; for characteristic $\ell$, see \cite{Ju-thesis}).
For example, in the case $G = GL_n$, the simple $k\SG_n$-modules are
denoted $D_\mu$, where $\mu$ runs over a subset of the partitions of
$n$, whose elements are called $\ell$-regular partitions, and it turns
out that the modular Springer correspondence is given, as in
characteristic zero, by the transposition of partitions
\cite[\S 6.4]{Ju-thesis}.

The fact that there is an inverse Fourier transform implies that
the map $\Psi_k$ is an injection. The fact that the Fourier transform
of the constant perverse sheaf is the sky-scraper sheaf concentrated
in zero implies that the pair corresponding to the trivial
representation consists of the trivial orbit and the trivial character.

Now, if $E$ is a $\KM W$-module, we can choose an $\OM$-lattice
$E_\OM$ stable by $W$. Then $\FM \otimes_\OM E_\OM$ is an $\FM
W$-module whose class in the Grothendieck group does not depend on the
choice of the lattice~\cite{Serre}. Thus we have well-defined multiplicities
$d^W_{E,F} = [\FM \otimes_\OM E_\OM : F]$ for $F$ an $\FM W$-module.
One can similarly define decomposition numbers
$d^\NC_{(x,\rho),(y,\sigma)}$ for $G$-equivariant perverse sheaves on
the nilpotent cone, where $(x,\rho) \in \PG_\KM$ and
$(y,\sigma)\in\PG_\FM$. Then we have
\[
d^W_{E,F} = d^\NC_{\Psi_\KM(E),\Psi_\FM(F)}
\]
which shows that the decomposition matrix of $W$ is a submatrix of the
decomposition matrix for $G$-equivariant perverse sheaves on the
nilpotent cone \cite{Ju-thesis}. It follows that the stalks of IC sheaves in
characteristic $\ell$ on the nilpotent singularities encode the
modular representation theory of Weyl groups.

\section{Perverse sheaves}

\subsection{Constructible sheaves}

Throughout, $k$ will denote a field or $\ZM$. All varieties will be
varieties over the complex numbers equipped with the classical
topology and all morphisms will be morphisms of varieties. Dimension
will always mean the complex dimension.

Let $X$ be a variety. We will denote by $\SS$ a decomposition
\begin{equation} \label{eq:dec}
X = \bigsqcup_{S \in \SS} S
\end{equation}
of $X$ into finitely many locally closed (in the Zariski topology)
connected smooth subvarieties. A sheaf of $k$-vector spaces $\FS$ on
$X$ will be called $\SS$-constructible if the restriction of $\FS$ to
each $S \in \SS$ is a local system (a sheaf of $k$-modules 
which is locally isomorphic to a constant sheaf with values in a finitely
generated $k$-module).  A sheaf $\FS$ is constructible if there exists an
$\SS$ as above making it $\SS$-constructible.

Let $D^b(X,k)$ denote the bounded derived category of sheaves of
$k$-vector spaces. Given $\KS \in D^b(X,k)$ we denote its cohomology
sheaves by $\HC^m(\KS)$. We denote by $D^b_c(X,k)$
(resp. $D^b_{\SS}(X,k)$) the full subcategory of $D^b(X,k)$ with
objects consisting of complexes $\KS \in D^b(X,k)$ such that
$\HC^m(\KS)$ is constructible (resp. $\SS$-constructible) for all
$m$. We have truncation functors $\tau_{\le i}$
and $\tau_{> i}$ on $D^b(X,k)$, $D^b_c(X,k)$ and $D^b_{\SS}(X,k)$. For
example, $\HC^m( \tau_{\le i} \KS)$ is isomorphic to $\HC^m( \KS)$ if
$m \le i$ and is 0 otherwise.

We have internal bifunctors $\RHOM$ and $\otimes_k^L$ on
$D^b_c(X,k)$, which are the derived functor of the usual bifunctors on
categories of sheaves. For any morphism $f : X \to Y$ we have functors:
\begin{equation*}
\xymatrix{ D^b_c(X,k) \ar@(ur,ul)[r]^{f_*, f_!} & D^b_c(Y,k)
  \ar@(dl,dr)[l]^{f^*, f^!} }
\end{equation*}
The functors $f_*$ and $f_!$, sometimes denoted $Rf_*$ and 
$Rf_!$, are the right derived functors of the direct image and 
direct image with compact support (both functors 
are left exact). The inverse image functor $f^*$ for sheaves is exact,
and passes trivially to the derived category. The pair $(f^*,f_*)$ is
adjoint. It turns out that the derived functor $f_!$ has a right
adjoint, namely $f^!$.

In the case where $Y$ is a point, for $\FS \in D^b_c(X,k)$ we have
\begin{align*}
f_*\FS &= \RGa(X,\FS), &
f^*\un k_\pt &= \un k_X, 
\\
f_!\FS &= \RGc(X,\FS),&
f^!\un k_\pt &= \DC_X,
\end{align*}
where $\DC_X$ is the dualizing sheaf of $X$. It allows one to define the
dualizing functor:
\begin{equation*}
\DM = \RHOM(-,\DC_X) : D^b(X,k)^\op \longto D^b(X,k).
\end{equation*}
Its square is isomorphic to the identity functor.
For example, if $X$ is smooth of dimension $d$ and $\LC$ is a local
system on $X$ then $\DM(\LC[d]) \simeq \LC^{\vee}[d]$, where
$\LC^{\vee}$ denotes the dual local system. In that case, we have
$\DC_X \simeq \uk_X[2d]$, and $\uk_X[d]$ is self-dual.

In general, we have isomorphisms $\DM f_* \simeq f_! \DM$ and $\DM f^* \simeq f^! \DM$. If $Y$
is a point, $X$ is smooth and $\LC$ is a local system then the first
isomorphism yields Poincaré duality between $H^*(X, \LC)$ and
$H^*_c(X, \LC^{\vee})$.

In what follows, we will fix a decomposition $\SS$ of 
$X$ as in (\ref{eq:dec}) and assume additionally that $\SS$ is a Whitney stratification.  
For $\SS$ to be a stratification, we require that the closure
of a stratum is a union of strata. The Whitney conditions,\footnote{
See \cite{Wh} p. 540, \cite{BoAl} A'Campo \S IV.1 p. 41, or
\cite{Li} \S 1--2. 
}
which we are not going to describe, will ensure that the functors
induced by inclusions of unions of strata, and the duality, preserve
the notion of $\SS$-constructibility. Any stratification of $X$ can be
refined into a Whitney stratification.

As an important special case, if $X$ is a $G$-variety with finitely
many orbits, where $G$ is a connected algebraic group, we can choose
for $\SS$ the set of $G$-orbits on $X$. In that case, we will be
interested in $G$-equivariant sheaves. This notion will be explained below.

\subsection{Perverse sheaves with coefficients in a field}
\label{subsec:perv field}

Throughout this section we assume that $k$ is a field. As in the
previous section we fix a variety $X$ with Whitney stratification
$\SS$. The category of perverse sheaves constructible with respect to
$\SS$, denoted $\Perv_\SS(X,k)$, consists of the full subcategory of
those objects $\FS \in D^b_{\SS}(X,k)$ such that:
\begin{enumerate}
\item\label{perv1} for all $S\in\SS$,
$i_S^* \FS$ is concentrated in degrees $\le - \dim S$,
\item\label{perv2} for all $S \in \SS$,
$i_S^! \FS$ is concentrated in degrees $\ge - \dim S$.
\end{enumerate}
Note that these two conditions are exchanged by $\DM$. It follows that
$\DM$ preserves $\Perv_{\SS}(X,k)$. We say that $\Perv_\SS(X,k)$ is
the heart of the $t$-structure
$(D^{\le 0}_\SS(X,k),D^{\ge 0}_\SS(X,k))$ where $D^{\le 0}_\SS(X,k)$,
resp. $D^{\ge 0}_\SS(X,k)$, is the full subcategory of $D^b_\SS(X,k)$
with objects satisfying condition \eqref{perv1}, resp. \eqref{perv2}.

For $\FS \in \Perv_{\SS}(X,k)$, we are interested in the stalks of the
cohomology sheaves of $\FS$. An induction shows that any perverse
sheaf $\FS$ only has non-trivial stalks in degrees $\ge - d$, where
$d = \dim X$. Hence, we see that $\FS$ is perverse if and only if the
cohomology sheaves of both $\FS$ and $\DM \FS$ along strata are of the
following form (see \cite{Ar}):
\begin{equation*}
\begin{array}{|c|c|c|c|c|c|c|}
\hline 
\mathrm{strata} & \dots & -d & -d+1 & \dots & -1 & 0 \\
 \hline 
S_{d} & 0 & * & 0 & 0 & 0 & 0 \\
 \hline 
S_{d-1} & 0 & * & * & 0 & 0 & 0 \\
 \hline 
\vdots & 0 & \vdots & \vdots & \ddots & 0 & 0 \\
 \hline 
S_1 & 0 & * & * & \dots & * & 0 \\
 \hline 
S_0 & 0 & * & * & \dots & * & * \\
 \hline
\end{array}
\end{equation*}
Here $S_m$ denotes the union of strata of dimension $m$ and $*$
denotes the possibility of a non-trivial cohomology sheaf. For
example, the first line tells us that the restriction of $\FS$ to an
open stratum is either zero or a local system concentrated degree
$-d$. The last line tells us that the stalks and costalks of $\FS$ at
0-dimensional strata can be non-trivial in degrees between $-d$ and
$0$.

The category $\Perv_{\SS}(X,k)$ is an abelian category and every object
has finite length. The exact sequences in $\Perv_{\SS}(X,k)$ are those
sequences
\begin{equation*}
\FS_1 \to \FS_2 \to \FS_3
\end{equation*}
which can be completed, via a map $\FS_3 \to \FS_1[1]$
(necessarily unique in this case), to distinguished
triangles in $D^b_{\SS}(X,k)$.

As in any abelian category of finite length, it is important to
understand the simple objects.  One has a bijection:
\begin{equation*}
\left \{
\begin{array}{c}
\text{ simple objects } \\ 
\text{ in } \Perv_{\SS}(X,k) \end{array}
\right \}
\elem{\sim}
\left\{
\begin{array}{c}
\text{ pairs $( S, \LC)$ where $S \in \SS$ }\\
\text{ and $\LC$ is an irreducible} \\
\text{ local system on $S$ } \end{array}
\right \}.
\end{equation*}
This bijection may be described as follows. Given a pair $(S, \LC)$
there exists a unique object $\ic(\overline{S}, \LC) \in
\Perv_{\SS}(X,k)$ such that
\begin{enumerate}
\item $i_S^* \ic(\overline{S}, \LC) \simeq \LC[d_S]$,
\item $\ic(\overline{S}, \LC)$ is supported on $\overline{S}$,
\end{enumerate}
and, for all strata $T \subset \overline{S}$ with $T \ne S$,
\begin{enumerate}
\item[(2)] $i_T^* \ic(\overline{S}, \LC)$ is concentrated in degrees
  $< -\dim T$,
\item[(3)] $i_T^! \ic(\overline{S}, \LC)$ is concentrated in degrees
  $> -\dim T$.
\end{enumerate}
The object $\ic(\overline{S}, \LC)$ is called the intersection
cohomology complex corresponding to $(S, \LC)$.
A different convention is to shift this complex by $-\dim S$, so that it
is concentrated in non-negative degrees. The normalization we use has
the advantage that the the intersection cohomology complexes are
perverse sheaves.

Note that we must have
$\DM \ic(\overline{S}, \LC) \simeq \ic(\overline{S}, \LC^{\vee})$.
This explains the existence of a Poincaré duality between
$\IH^*(X, \LC)$ and $\IH_c^*(X, \LC^{\vee})$.

It is useful to note the special form of the restrictions to the
strata of an intersection
cohomology complex $\ic(\overline{S}, \LC)$, depicted as before (note
the zeroes on the diagonal). Let $d_S$ denote the dimension of $S$, so
that $S \subset S_{d_S}$.

\begin{equation*}
\begin{array}{|c|c|c|c|c|c|c|c|c|}
 \hline
\mathrm{strata} & - d & \dots & - d_S - 1 & -d_S & - d_S + 1 & \dots & -1 & 0\\
 \hline
S_{d} & {\bf 0} & 0 & 0 & 0 & 0 & 0 & 0 & 0\\
 \hline
\vdots & 0 & {\bf 0} & 0 &0 & 0 & 0 & 0 & 0\\
 \hline
S_{d_S + 1} & 0 & 0 & {\bf 0} & 0 & 0 & 0 & 0 & 0 \\
 \hline
S_{d_S} & 0 & 0 & 0 & {\bf i_{S*}\LC} & 0 & 0 & 0 & 0 \\
 \hline
S_{d_S - 1} & 0 & 0 & 0 & * & {\bf 0} & 0 & 0 & 0 \\
 \hline
\vdots & 0 & 0 & 0 & * & * & {\bf 0} & 0 & 0 \\
 \hline
S_1 & 0 & 0 & 0 & * & * & * & {\bf 0} & 0\\
 \hline
S_0 & 0 & 0 & 0 & * & * & * & * & {\bf 0}\\
 \hline
\end{array}
\end{equation*}

\bigskip

Let us make a small digression about equivariance. We follow
\cite[\S 0]{ICC}. Assume $X$ is endowed with an action of a connected
algebraic group $G$. Let $a:G\times X \to X$ and
$\pr:G\times X \to X$ denote the action morphism and the second
projection. A perverse sheaf $\FS$ is called $G$-equivariant if
there is an isomorphism $\a:a^*\FS \xra{\sim} \pr^* \FS$. This isomorphism is
unique if we impose the condition that the induced isomorphism
$i^*a^*\FS \xra{\sim} i^*\pi^*\FS$ is the identity of $\FS$, where
the morphism $i : X \to G \times X$, $x\mapsto (1,x)$ is a section of
$a$ and $\pi$, so that $i^*a^*\FS =\FS$ and $i^*\pi^*\FS = \FS$.
This follows from \cite[Prop. 4.2.5]{BBD} (Deligne). Then $\a$
satisfies the usual associativity condition. Note that this definition
of $G$-equivariance for perverse sheaves does not work for arbitrary
complexes.\footnote{
One should instead consider the equivariant 
derived category as defined in \cite{BLu}.}
It also works, though, for usual sheaves (complexes
concentrated in one degree), and in particular for local systems.
The $G$-equivariant local systems on an orbit $S$ correspond
bijectively to the finite dimensional representations of the finite
group of components of the isotropy group of a point $x$ in $S$.
We are particularly interested in the case where $G$ has finitely many
orbits in $X$, and we take the stratification of $X$ into its
$G$-orbits. Then the simple $G$-equivariant perverse sheaves
correspond to pairs consisting of an orbit and an irreducible
$G$-equivariant local system on that orbit.

We now recall Deligne's construction of intersection cohomology
complexes. As above let $S_m$ denote the union of strata of dimension
$m$ and denote by $X_m$ the union of all strata of dimension greater
than or equal to $m$. We have a sequence of inclusions:
\[
X_d \stackrel{j_{d-1}}{\hookrightarrow} X_{d-1}
\stackrel{j_{d-2}}{\hookrightarrow} X_{d-2} \hookrightarrow \cdots
\hookrightarrow X_1 \stackrel{j_0}{\hookrightarrow} X_0 = X.
\]

Now let $\LC$ be an irreducible local system on $S \in \SS$. We still denote
by $\LC$ its extension by zero to $X_{d_S}$. One has an
isomorphism
\[
\ic(\overline{S}, \LC) \simeq
(\tau_{\le -1} \circ j_{0*}) \circ
(\tau_{\le -2} \circ j_{1*}) \circ
\dots \circ
(\tau_{\le - d_S} \circ j_{d_S - 1*}) (\LC[d_S]).
\]

This allows the calculation of $\ic(X, \LC)$ inductively on the strata. We
will see examples of this construction below. However, the $j_*$
functors are not easy to compute explicitly in general.

In characteristic zero, the decomposition theorem provides
a much more powerful means of calculating the
stalks of intersection cohomology complexes.
Given a Laurent polynomial $P = \sum a_i v^i
\in \NM[v,v^{-1}]$ and $\KS \in D^b(X,k)$, we define
\begin{equation*}
P \cdot \KS = \bigoplus \KS[i]^{\oplus a_i}.
\end{equation*}
We call a complex $\KS \in D^b_{\SS}(X,k)$ semi-simple if one has an
isomorphism
\begin{equation*}
\KS \simeq \bigoplus P_{S, \LC_S} \cdot \ic(\overline{S}, \LC_S)
\end{equation*}
for some $P_{S, \LC_S} \in \NM[v,v^{-1}]$, where the sum is over all pairs $(S,
\LC)$ consisting of an irreducible local system on a stratum $S$.
In other words, an object $\KS \in D^b_{\SS}(X,k)$ is semi-simple if it
is isomorphic to a direct sum of shifts of irreducible intersection
cohomology complexes.
One version of the decomposition theorem is the following:

\begin{thm}[\cite{BBD}, \cite{Saito}, \cite{dCM}]
Let $\pi : Y \to X$ be a proper map from a smooth variety $Y$. If $k$
is of characteristic zero, then $\pi_* \uk_Y[\dim Y]$ is semi-simple.
\end{thm}

We will see examples below of the failure of the decomposition theorem
in positive characteristic.

Some conditions on the dimensions of the fibers of a resolution,
however, have consequences which hold for arbitrary coefficients.
Let $\pi:Y \to X$ be a proper morphism between $n$-dimensional
irreducible varieties, and assume $X$ is endowed with a stratification
$\SS$ such that $\pi$ is a weakly stratified mapping, that is, for each
stratum $S$ in $\SS$, the restriction of $\pi$ to $\pi^{-1}(S)$ is a
topological fibration with base $S$ and fibre $F_S$. Then $\pi$ is said
to be semi-small if, for all $S$, we have
$\dim F_S \leq \frac 1 2 \codim_X(S)$. A stratum $S$ is relevant for
$\pi$ if equality holds. We say that $\pi$ is small if the only relevant
stratum is the dense one. The following proposition is well-known
\cite{BM1,BM2,GM2}:

\begin{prop}
Let $\pi : Y \to X$ be a proper morphism as above with $Y$ smooth.
\begin{enumerate}
\item
If $\pi$ is semi-small then
$\pi_* \uk_Y[\dim Y]$ is a perverse sheaf.
\item
If $\pi$ is small then $\pi_* \uk_Y[\dim Y]$ is an
intersection cohomology complex.
\end{enumerate}
\end{prop}

\subsection{Perverse sheaves over the integers}

In this section, we will give a flavour of the subtleties that occur
when we take integer coefficients. See \cite[\S 3.3]{BBD}, and
\cite{decperv} for a more detailed study.

Let us first consider the bounded derived category of constructible
sheaves of $\ZM$-modules on the point $\pt = \Spec \CM$. These are
just complexes of $\ZM$-modules with finitely many non-zero cohomology groups,
all of which are finitely generated over $\ZM$. The duality functor is
$\DM = \RHom(-,\ZM)$. The perverse sheaves on $\pt$ are just finitely
generated $\ZM$-modules placed in degree zero. They form the heart of
the usual $t$-structure, corresponding to the perversity $p$. We will
see that this $t$-structure is not preserved by the duality.

Since $\ZM$ is hereditary (of global dimension $1$), any object of
$D^b_c(\pt,\ZM)$ is isomorphic to the direct sum of its shifted
cohomology objects. So the indecomposable objects in $D^b_c(\pt,\ZM)$
are concentrated in one degree, and they are isomorphic, up to shift,
either to $\ZM$ or to $\ZM / \ell^a$ for some prime number $\ell$, and
some positive integer $a$.

First consider the indecomposable $\ZM$. It is a free, hence projective
$\ZM$-module, and thus we have
$\DM(\ZM) = \RHom(\ZM,\ZM) = \Hom(\ZM,\ZM) = \ZM$.
Here there is no problem. The dual of a torsion-free module remains in
the heart of the natural $t$-structure.

Now consider the case of $\ZM/\ell^a$. Here we cannot apply the
functor $\Hom(-,\ZM)$ directly to the module $\ZM/\ell^a$. First we
have to replace $\ZM/\ell^a$ by a projective resolution, as $\ZM
\xra{\ell^a} \ZM$, with the first $\ZM$ in degree $-1$. Now we can
apply the functor $\Hom(-,\ZM)$, and we get the complex $\ZM
\xra{\ell^a} \ZM$ with the last $\ZM$ in degree $1$. So the dual of
$\ZM/\ell^a$ is isomorphic to $\ZM/\ell^a [-1]$. This is another way
to say that $\Hom(\ZM/\ell^a, \ZM) = 0$,
$\Ext^1(\ZM/\ell^a, \ZM) \simeq \ZM/\ell^a$,
and $\Ext^i(\ZM/\ell^a, \ZM) = 0$ for $i > 1$.

This shows that the usual $t$-structure (for the perversity $p$) is
not stable by the duality.  We see that the problem comes from the
torsion. The duality exchanges a free module in degree $n$ with a free
module in degree $-n$, but it exchanges a torsion module in degree $n$
with a torsion module in degree $1 - n$. The duality exchanges the
usual $t$-structure $(D^{\le 0},D^{\ge 0})$ on
$D^b_c(\pt,\ZM)$ with another $t$-structure
$(D^{\le 0^+},D^{\ge 0^+})$, defined by:
\begin{gather*}
K \in D^{\le 0^+} \Iff
H^1(K) \text{ is torsion and } H^i(K) = 0 \text{ for i > 1}\\
K \in D^{\ge 0^+} \Iff
H^0(K) \text{ is torsion-free and } H^i(K) = 0 \text{ for i < 0}
\end{gather*}

Remember that we constructed the $t$-structure for the perversity $p$
by taking on each stratum $S$ the usual $t$-structure shifted by
$\dim S$, and gluing them together. Over the integers, we can
either do the same, or take on each stratum the dual of the usual
$t$-structure, shifted by $\dim S$ (here we consider torsion versus
torsion-free local systems on $S$), and then again glue them together.
In the second case, we obtain the following $t$-structure,
corresponding to the perversity $p_+$:
\begin{gather*}
\FS\in D^{\le 0^+}_\SS(X,\ZM) \Iff
\begin{cases}
\HC^m(i_S^* \FS) \text{ is zero for } m > - \dim S + 1\\
\text{ and is torsion for } m = -\dim S + 1
\end{cases}
\\
\FS\in D^{\ge 0^+}_\SS(X,\ZM) \Iff
\begin{cases}
\HC^m(i_S^! \FS)= 0 \text{ is zero for } m < - \dim S\\
\text{ and torsion-free for } m = -\dim S
\end{cases}
\end{gather*}

We will denote the heart of the classical $t$-structure
by $\Perv_{\SS}(X,\ZM)$, and the heart of this new $t$-structure by
$\Perv^+_{\SS}(X,\ZM)$.

The abelian category $\Perv_{\SS}(X,\ZM)$ is Noetherian but not Artinian
(again, this is already the case for $X = pt$). However, given any
local system $\LC$ on a stratum $S$ one still has a unique extension
$\ic(\overline{S}, \LC)$ satisfying the same conditions as for a
field. Keeping the notation of the previous section, this may be
defined by:
\[
\ic(\overline{S}, \LC) :=
(\tau_{\le -1} \circ j_{0*}) \circ
(\tau_{\le -2} \circ j_{1*}) \circ
\dots \circ
(\tau_{\le - d_S} \circ j_{d_S - 1*}) (\LC[d_S]).
\]

To obtain the dual of an intersection cohomology complex
$\ic(\overline{S}, \LC)$ in $\Perv^+_{\SS}(X, k)$ one needs to
consider a variant of the truncation functors on $D^b(X, \ZM)$, which
we denote $\tau^+_{\le i}$ and $\tau^+_{> i}$. If $\KS \in D^b(X,\ZM)$
then $\HC^m(\tau^+_{\le i} \KS)$ is isomorphic to $\HC^m(\KS)$ for $m
\le i$, to the torsion submodule of $\HC^m(\KS)$ for $m = i + 1$, and
is zero otherwise. One defines
\begin{equation*}
\ic^+(\overline{S}, \LC) :=
(\tau^+_{\le -1} \circ j_{0*}) \circ
(\tau^+_{\le -2} \circ j_{1*}) \circ
\dots \circ 
(\tau^+_{\le - d_S} \circ j_{d_S - 1*}) (\LC[d_S]).
\end{equation*}
As $\ic(\overline{S}, \LC)$, the complex $\ic^+(\overline{S}, \LC)$ may be
characterized in terms of the stalks of $i_S^*$ and $i_S^!$.
As one might expect, the complexes
$\ic(\overline{S}, \LC)$ and $\ic^+(\overline{S},\LC^{\vee})$
are exchanged by the duality.

\subsection{First example: the nilpotent cone of $\sl_2$}
\label{subsec:sl_2}

Let ${\sl}_2$ be the Lie algebra of $2\times 2$ traceless matrices over
$\CM$ and let $\NC \subset {\sl}_2$ be its nilpotent cone. It is
isomorphic to a quadratic cone inside affine 3-space:
\begin{equation*}
\NC =
\left \{
\begin{pmatrix}
x & y \\
z & -x 
\end{pmatrix}
\middle| \;
x^2 + yz = 0 
\right\}
\subset
{\sl}_2 \simeq
{\AM^3}
\end{equation*}
Note also that $\NC$ is isomorphic to the quotient of a two
dimensional vector space $V = \Spec \CM[u,v]$ by the scalar action of
$\{ \pm 1 \}$. If we
choose coordinates $(u,v)$ on $V$ then an isomorphism is given by
\begin{align*}
V / \{ \pm 1 \} & \longto \NC \\
\pm (u,v) &\longmapsto \left (\begin{array}{cc} uv
  & -v^2 \\ u^2 & -uv \end{array} \right ).
\end{align*}
The conjugation action of $SL_2(\CM)$ on $\NC$ has two orbits,
$\OC_\reg$ and $\{ 0 \}$, and we let $\SS$ denote the stratification of
$\NC$ into these two orbits. We will be interested in calculating
$\ic(\NC, k)$ for $k = \QM$, $\ZM$ and $\FM_p$ (of course $\ic(\{ 0 \}
, k)$ is always a skyscraper sheaf on $\{0 \}$ in degree 0).

We will first examine Springer's resolution
\begin{equation*}
T^* \PM^1 \to \NC.
\end{equation*}
Concretely, we may identify $T^* \PM^1$ with pairs $( \ell, x)$ where
$\ell \in \PM^1$ is a line containing the image of $x \in \NC$. The
map $\pi$ is then obtained by forgetting $\ell$; it is clearly an
isomorphism over $\OC_\reg$ and has fibre $\PM^1$ over $\{0 \}$. Hence,
for any $k$ the stalks of $\pi _* \uk_{T^* \PM^1}[2]$ are given by
\begin{equation*}
\begin{array}{|c|c|c|c|}
\hline & -2 & -1 & 0 \\ \hline \OC_\reg & k & 0 & 0\\ \hline \{ 0 \} &
k & 0 & k \\ \hline
\end{array}
\end{equation*}
If $k$ is of characteristic 0, then we know by the decomposition
theorem that $\pi _* \uk_{T^* \PM^1}[2]$ is semi-simple and hence
\begin{equation*}
\pi _* \uk_{T^*\PM^1}[2] \simeq \ic(\NC,k) \oplus \ic(\{0 \}, k).
\end{equation*}
It follows that $\ic(\NC, k)$ is isomorphic to $\uk_{\NC}[2]$.

To handle the case $k = \FM_p$ requires more care. Recall that the
Deligne construction tells us that
\begin{equation*}
\ic(\NC,k) \simeq \tau_{\le -1} \circ j_* (\uk_{\OC_\reg} [2]),
\end{equation*}
where $j : \OC_\reg \injto \NC$ is the open immersion.
Since this complex is an extension of $\uk_{\OC_\reg} [2]$, the only
stalk that we have to compute is the stalk at zero. Let us compute
$(j_* \uk_{\OC_\reg})_0$. Then it will be a trivial matter to shift
and truncate.

We have
\[
(j_* \uk_{\OC_\reg})_0 \simeq
\lim_{V \ni 0} \RGa(V \setminus \{ 0 \}, k)
\]
where $V$ runs over the open neighbourhoods of $0$ in $\NC$.

We can replace this limit by a limit over a basis of neighbourhoods of
zero, as for example $V_n$, $n\ge 1$, the intersection of $\NC$ with
the open ball of radius $1/n$ centered at $0$ in $\CM^3$.
But $\NC$ is a cone: it is stable by multiplication by a scalar in
$\CM^*$, and in particular in $\RM_{>0}$. We can use this to see that
all the $V_n \setminus \{0\}$ are homeomorphic: they are actually
homeomorphic to $\NC \setminus \{0\} = \OC_\reg$ itself.
Thus we have
\begin{equation*}
(j_* \uk_{\OC_\reg})_0 \simeq 
\lim_n \RGa(V_n \setminus \{ 0 \}, k) \simeq
\RGa(\OC_\reg, k).
\end{equation*}
(this argument applies for any cone).
However, as we observed above, $\NC \simeq V/\{\pm 1\}$ and hence
$\OC_\reg \simeq (V \setminus \{ 0 \} )/ \{\pm 1\}
= (\CM^2\setminus \{0\}) / \{\pm 1\}$, 
which is homotopic to $S^3 / \{\pm 1\} = \RM\PM^3$.
Thus
\begin{equation*}
(j_* \uk_{\OC_\reg})_0 \simeq 
\RGa(\RM\PM^3, k) \simeq
(k \elem{0} k \elem{2} k \elem{0} k),
\end{equation*}
the latter complex being concentrated in degrees between $0$ and $3$.

Hence, if $k$ is a field of characteristic $p$, the stalks of
$j_*\uk_{\OC_\reg}[2]$ are given by:
\begin{equation*}
\begin{array}{|c|c|c|c|c|}
\hline & -2 & -1 & 0 & 1 \\
 \hline 
\OC_\reg & k & 0 & 0 & 0\\
 \hline
\{ 0 \} & k & (k)_2 & (k)_2 & k \\
 \hline
\end{array}
\end{equation*}
where $(k)_2$ means $k$ if $p = 2$, and $0$ otherwise.
We obtain the stalks of $\ic(\NC, k)$ by truncating:
\begin{equation*}
\begin{array}{|c|c|c|c|}
\hline & -2 & -1 & 0 \\
\hline \OC_\reg & k & 0 & 0\\ 
\hline \{ 0 \} & k & (k)_2 & 0\\
\hline
\end{array}
\end{equation*}
In fact, the decomposition theorem holds here if and only if $p \ne 2$.

One may calculate the stalks of $\ic(\NC, \ZM)$ and 
$\ic^+(\NC, \ZM)$ and one obtains:
\begin{equation*}
\begin{array}{|c|c|c|c|}
\hline & -2 & -1 & 0 \\ \hline \OC_\reg & \ZM & 0 & 0 \\ \hline \{ 0 \}
& \ZM & 0 & 0 \\ \hline
\end{array}
\qquad
\begin{array}{|c|c|c|c|}
\hline & -2 & -1 & 0 \\ \hline \OC_\reg & \ZM & 0 & 0 \\ \hline \{ 0 \}
& \ZM & 0 & \ZM/2\ZM \\ \hline
\end{array}
\end{equation*}

>From this information, one calculate the decomposition numbers for
$GL_2$-equivariant perverse sheaves on its nilpotent cone \cite{decperv}, and we get:
\[
\begin{array}{lc}
& (1^2) \ (2)\\
\begin{array}{c}
(1^2)\\
(2)
\end{array}
&
\begin{pmatrix}
{\bf 1}&0\\
{\bf 1}&1
\end{pmatrix}
\end{array}
\]
from which we can extract the decomposition matrix for $\SG_2$ in
characteristic $2$:
\[
\begin{array}{lc}
& D_{(2)}\\
\begin{array}{c}
S_{(2)}\\
S_{(1^2)}
\end{array}
&
\begin{pmatrix}
{\bf 1}\\
{\bf 1}
\end{pmatrix}
\end{array}
\]

\section{Some stalks in nilpotent cones}

In this section we give more examples of calculations of stalks of
intersection cohomology sheaves on nilpotent cones, motivated by the
modular Springer correspondence. Actually, in type $A$ the nilpotent
singularities also occur in the affine Grassmannian \cite{lu}, so the
geometric Satake theorem is another motivation.

First we recall a parabolic generalization of Springer's resolution,
providing a resolution of Richardson orbit closures. Then we deal with
the minimal nilpotent orbit closure in ${\sl}_n$ and ${\sp}_{2n}$, and
with the singularity of the nilpotent cone of ${\sl}_n$ at the
subregular orbit, by a direct approach. Note that minimal and simple
singularities have been dealt with in all types in \cite{cohmin,
  decperv}, but in the special cases we treat here, the calculation
can be done quickly. Finally, we give more computations of stalks of
intersection cohomology complexes in characteristic $p$, which might be
new: we give all the stalks in the nilpotent cone of ${\sl}_3$ for
$p\neq 3$, and all the stalks in the subvariety of the nilpotent cone
of ${\sl}_4$ consisting of the matrices which square to zero.

All our calculations will be completed for varieties over $\CM$ in the metric
topology, however all the results and most proofs can be translated
into the étale situation: we mostly use basic facts about the
cohomology of projective spaces and flag varieties, Gysin sequences
etc. which have direct translations in \'etale cohomology. We will
always outline how a calculation can be performed in the \'etale
topology when such a direct translation is not possible.

Throughout we will use the following notation (already used in the
last section). If $k$ is a field of characteristic $p$ and $n$ is an
integer then
\begin{equation*}
(k)_n = 
\begin{cases}
k & \text{if $p$ divides $n$,} \\
0 & \text{otherwise.}
\end{cases}
\end{equation*}

\subsection{Semi-small resolutions of Richardson orbit closures}
\label{subsec:richardson}

Remember the Springer resolution:
\[
\begin{array}{rcccc}
\pi_\NC : T^*(G/B) &=& G \times^B \ug &\longto& \NC = \ov \OC_\reg\\
&&g *_B x &\longmapsto& (\Ad g)(x)
\end{array}
\]
We will see a generalization of this resolution, where we replace the
Borel subgroup $B$ by a parabolic subgroup $P$, with unipotent radical
$U_P$. We denote by $\ug_P$ the Lie algebra of $U_P$. We can naturally
define the proper morphism:
\[
\begin{array}{rcccc}
\pi_P : T^*(G/P) &=& G \times^P \ug_P &\longto& \NC\\
&&g *_P x &\longmapsto& (\Ad g)(x)
\end{array}
\]
but what is the image of $\pi_P$~?

Since $\ug_P$ is irreducible and there are only finitely many
nilpotent orbits, there is a unique nilpotent orbit $\OC$ such that
$\OC \cap \ug_P$ is dense in $\ug_P$. This orbit is called the
Richardson orbit associated to $P$. One can see that the image of
$\pi_P$ is the closure of $\OC$.  The morphism $\pi_P$ induces a
semi-small resolution of $\ov \OC$, which is used in \cite{BM1}.
In type $A$, all nilpotent orbits are Richardson, but in general this
is not so. The regular nilpotent orbit, though, is always Richardson:
it is the one associated to $B$.

Let us describe the situation in the case $G = SL_n$.  If $\l =
(\l_1,\dots,\l_s)$ is a partition of $n$,
let $\OC_\l$ denote the nilpotent orbit consisting of the
nilpotent matrices whose Jordan normal form $x_\l$ has Jordan blocks
of sizes given by the parts of $\l$. We denote by $P_\l$ the parabolic
subgroup of $SL_n$ stabilizing the standard partial flag of shape $\l$:
\[
F^\l_\bullet := (0 \subset F^\l_1 \subset \dots \subset F^\l_s = \CM^n)
\]
where $F^\l_i$ is spanned by the $\l_1 + \dots + \l_i$ first elements
of the canonical basis of $\CM^n$. The partial flag variety $G/P_\l$
can be interpreted as the variety $\FC_\l$ of all partial flags of shape $\l$,
that is, sequences of subspaces
$F_\bullet = (0 \subset F_1 \subset \dots \subset F_s = \CM^n)$
with $\dim F_i = \l_1 + \dots + \l_i$. We have
\[
\begin{array}{rcl}
G \times^P \ug_{P_\l} &\elem{\sim}&
\{ (x,F_\bullet) \in \NC \times \FC_\l \mid x(F_i) \subset F_{i-1}\}\\
g *_{P_\l} x & \longmapsto & (\Ad g (x), g (F^\l_\bullet))
\end{array}
\]

For $x\in\NC$, we have
\[
x \in \OC_\l \Iff \forall i,\ \dim \Ker x^i = \l'_1 + \dots + \l'_i\\
\]
and
\[
x \in \ov \OC_\l \Iff \forall i,\ \dim \Ker x^i \geq \l'_1 + \dots + \l'_i\\
\]
where $\l'$ is the partition conjugate to $\l$.
Let us note that we have
\[
\OC_\mu \subset \ov \OC_\l
\Iff \forall i,\ \mu'_1 + \dots + \mu'_i \geq \l'_1 + \dots + \l'_i
\Iff \mu' \geq \l'
\Iff \mu \leq \l
\]
where $\leq$ is the usual dominance order on partitions.

Now assume $x$ is in the image of $\pi_{P_{\l'}}$. Then there is a flag
$F_\bullet$ of type $\l'$ such that $x(F_i) \subset F_{i-1}$ for all $i$. In
particular, we have $F_i \subset \Ker x^i$ for all $i$, and hence
$\dim \Ker x^i \geq \l'_1 + \dots \l'_i$ for all $i$. Thus
$x \in \ov\OC_\l$. Consequently, the image of $\pi_{P_{\l'}}$ is
included in $\ov\OC_\l$.

If $x \in \OC_\l$, then there is a unique flag $F_\bullet$ of type
$\l'$ such that $x(F_i) \subset F_{i-1}$, namely
$F_i = \Ker x^i$. Thus the image of $\pi_{P_{\l'}}$ contains $\OC_\l$,
and hence $\ov\OC_\l$, as it is a proper morphism.

Thus the image of $\pi_{P_{\l'}}$ is equal to $\ov\OC_\l$, and
$\pi_{P_{\l'}}$ is an isomorphism over $\OC_\l$. Since it is proper,
it is a resolution of singularities.

Note that this gives, in principle, a method to compute all the IC
stalks with $\QM$ coefficients of closures of nilpotent orbits in
${\sl}_n$, by induction, using the decomposition theorem: the direct image
$\pi_{P_{\l'}*} \un \QM [\dim \OC_\l]$ decomposes as a direct sum of
$\ic(\ov\OC_\l,k)$ and some copies of IC sheaves for lower strata,
which we know by induction. The stalks of the direct image are given
by the cohomology of the fibers. One finds the stalks of the IC sheaf
of $\ov\OC_\l$ by removing the stalks of the other summands.
It is a nice exercise to do that for small ranks. We will see some
examples below. Note, however, that all stalks of $G$-equivariant IC
complexes on nilpotent cones are known in characteristic zero, as
there is an algorithm to compute them \cite[V]{CS}, as soon as one has
determined the generalized Springer correspondence defined in
\cite{ICC}, which has also been done, by work of several authors.
In the case of $GL_n$, the answer (which is in terms of Kostka
polynomials) has been known since \cite{lu}.

With $\ZM$ or $\FM_p$ coefficients, however, one cannot use the
decomposition theorem, and the calculations are much more difficult.
To find the answer in general is a very deep and important problem; in
particular, such information would be sufficient to determine the
decomposition matrices of the symmetric groups, which is a central
problem in the modular representation theory of finite groups. We will
see some examples below where the calculations can be done.

\subsection{Minimal class in ${\sl}_n$}
\label{subsec:sl_n}

The goal of this paragraph will be to generalize the calculation \ref{subsec:sl_2}
of the stalk at the origin of the IC sheaf on the nilpotent cone for
${\sl_2}$ to that of the closure of the minimal non-trivial nilpotent
orbit in ${\sl_n}$. The minimal orbit $\OC_\mini = \OC_{(2,1^{n-2})}$ is the set of
nilpotent matrices with $1$-dimensional image and $(n-1)$-dimensional
kernel. Let us apply the considerations of the last subsection. Here
$\l = (2,1^{n-2})$, and $\l' = (n-1,1)$. The parabolic subgroup
$P = P_{(n-1,1)}$ is the stabilizer of a hyperplane. The partial flag
variety $G/P$ is identified with the projective space $\PM^{n-1}$ of
hyperplanes in $\CM^n$. We get a resolution of 
$\ov \OC_\mini$ by taking pairs $(x,H)$ consisting of a nilpotent
element $x$ (necessarily in $\ov \OC_\mini$) and a hyperplane $H$ contained in the
kernel of $x$. (Dually, we could have considered pairs
$(x,\ell)$ where $x$ is nilpotent and $\ell$ is a line such
that $\mathrm{Im}(x)\subset \ell$.)

Thus we have a proper morphism:
\begin{equation*}
T^* \PM^{n-1} \to \ov\OC_\mini = \OC_\mini \cup \{0\}
\end{equation*}
which is an isomorphism over $\OC_\mini$ and has
fibre $\PM^{n-1}$ over $\{0\}$. Hence, for any $k$ the stalks of
$\pi_* \uk_{T^* \PM^{n-1}}[2n-2]$ are given by:
\begin{equation*}
\begin{array}{|c|c|c|c|c|c|c|}
 \hline 
& - 2n + 2& - 2n + 1 & - 2n + 2 & \ldots & -1 & 0 \\
 \hline 
\OC_\mini & k & 0 & 0& \ldots & 0 & 0 \\
 \hline 
\{ 0 \} & k & 0 & k & \ldots & 0 & k\\
 \hline
\end{array}
\end{equation*}

If $k$ is of characteristic 0, then we know by the decomposition
theorem that $\pi _* \uk_{T^* \PM^{n-1}} [2n-2]$ is semi-simple and hence
\begin{equation*}
\pi _* \uk_{T^*\PM^{n-1}} [2n-2] \simeq \ic(\ov\OC_\mini,k) \oplus \ic(\{0 \}, k).
\end{equation*}
It follows that the cohomology of $\ic(\ov\OC_\mini, k)_0$ is
isomorphic to $k$ in even degrees between $-2n + 2$ and $-2$, and zero
otherwise. 

As in the ${\sl}_2$ case, to handle the cases $k = \FM_p$ or $\ZM$, we
use Deligne's construction:
\begin{equation*}
\ic(\ov\OC_\mini,k) \simeq \tau_{\le -1} \circ j_* (\uk_{\OC_\mini}[2n-2])
\end{equation*}
where $j$ is the open immersion $\OC_\mini \injto \ov\OC_\mini$.
Again, $\ov\OC_\mini$ is a cone and thus we have:
\begin{equation*}
(j_* \uk_{\OC_\mini})_0 \simeq \RGa(\OC_\mini, k).
\end{equation*}

Now we have observed that $\OC_\mini$ is isomorphic to the complement
of the zero section to the cotangent bundle of $\PM^{n-1}$.  The
cohomology can thus be read off the Gysin sequence:
\[
\cdots \longto
H^{i-2n+2} \elem{e} 
H^i \longto 
H^i(\OC_\mini) \longto
 H^{i-2n+3} \elem{e} 
H^{i+1} \longto
\cdots
\]
where $H^i := H^i(\PM^{n-1},\ZM)$
and $e$ is the Euler class of the cotangent
bundle. But this is simply the Euler characteristic
$n = \chi(\PM^{n-1})$ times a generator of
$H^{2n-2}(\PM^{n-1},\ZM) \simeq \ZM$.
We deduce that $H^i(\OC_\mini,\ZM)$ is isomorphic to $\ZM$ for
$i = 0,\ 2,\ \dots,\ 2n-4$, to $\ZM/n$ for $i = 2n - 2$,
and then again to $\ZM$ for $i = 2n - 1,\ 2n + 1,\ \dots,\ 4n - 5$,
and zero otherwise:
\begin{equation*}
\begin{array}{|c|c|c|c|c|c|c|c|c|c|c|}
\hline
0 & 1 & 2 & \cdots & 2n-3 & 2n-2 & \cdots & 4n-3 & 4n-4 & 4n-5 \\ 
\hline 
\ZM & 0 & \ZM & \cdots & 0 & \ZM/n & \cdots & \ZM & 0 & \ZM \\
 \hline
\end{array}
\end{equation*}

We have $\RGa(\OC_\mini,k) = k \otimes^L_\ZM \RGa(\OC_\mini,\ZM)$.
Thus each copy of $\ZM$ is replaced by $k$.  For $k = \QM$ or $\FM_p$
with $p \nmid n$, the torsion group $\ZM/n$ in degree $2n-2$ is
killed. For $k = \FM_p$ with $p\mid n$, $\ZM/n$ is replaced by two
copies of $\FM_p$, one in degree $2n-3$, one in degree $2n-2$.

We obtain the stalks of $\ic(\ov\OC_\mini, k)$ by shifting and
truncating. For $k = \FM_p$, we get:
\begin{equation*}
\begin{array}{|c|c|c|c|c|c|c|c|}
\hline & -2n+2 & -2n+3 & -2n+4 & \cdots & -2 & -1 & 0 \\
 \hline 
\OC_\mini & k & 0 & 0 & \cdots & 0 & 0 & 0 \\
 \hline 
\{0\} & k & 0 & k & \cdots & k & (k)_n & 0 \\
 \hline
\end{array}
\end{equation*}

For $k = \ZM$, we get:
\begin{equation*}
\begin{array}{|c|c|c|c|c|c|c|c|}
\hline & -2n+2 & -2n+3 & -2n+4 & \cdots & -2 & -1 & 0 \\
 \hline 
\OC_\mini & \ZM & 0 & 0 & \cdots & 0 & 0 & 0 \\
 \hline 
\{0\} & \ZM & 0 & \ZM & \cdots & \ZM & 0 & 0 \\
 \hline
\end{array}
\end{equation*}
For $\ic^+(\OC_\mini,k)$, one adds a copy of $\ZM/n$ in degree $0$ at $\{0\}$.

The above calculations give a geometric proof that the natural
representation of $\SG_n$ (which corresponds to the minimal orbit)
remains irreducible modulo $\ell$ if and only if $\ell \nmid n$, while
if $\ell \mid n$ its modular reduction involves the trivial
representation once.

\subsection{Minimal class in $\sp_{2n}$}
\label{subsec:sp_2n}

In this section, we will treat the case of the minimal class
$\OC_\mini$ in $\gg = \sp_{2n}$. We view $G = Sp_{2n}$ as the subgroup
of $GL_{2n}$ stabilizing some symplectic form on $\CM^{2n}$ given by a
matrix $Q$ with respect to the canonical basis.

A matrix $M \in GL_{2n}$ will be in $Sp_{2n}$ if and only if ${}^t M Q
M = Q$.  Then $\gg$ can be identified with the Lie algebra of the
matrices $H \in \gl_{2n}$ such that the following identity holds
in $GL_{2n}(\CM[\varepsilon])$: 
\[
(1 + \e {}^t H) Q (1 + \e H) = Q
\]
where $\CM[\varepsilon] = \CM[X]/(X^2)$ and
$\varepsilon$ is the image of $X$ in this quotient.  This is
equivalent to
\begin{equation}
\label{sp_2n}
{}^tHQ + QH = 0
\end{equation}
in $\gl_{2n}$.  Now, the minimal class $\OC_\mini$ consists of those
matrices in $\sp_{2n}$ which are nilpotent with Jordan type
$(2,1^{2n-2})$. They are also characterized in $\sp_{2n}$ by the fact
that they are of rank one (all matrices in $\sp_{2n}$ have zero
trace).

A matrix $H \in \gl_{2n}$ of rank one is of the form $H = u{}^tv$,
where $u$ and $v$ are non-zero vectors in $\CM^{2n}$.  Moreover, $u$
and $v$ are uniquely determined up to multiplying $u$ by some non-zero
scalar $\lambda$ and dividing $v$ by the same scalar $\lambda$.  Now
suppose $H$ is in $\OC_\mini$. Then \eqref{sp_2n} writes
\[
v{}^tu Q + Q u{}^tv = 0
\]
that is,
\[
(Qu){}^tv = v{}^t(Qu)
\]
using the fact that $Q$ is anti-symmetric.  Since $Q$ is
non-degenerate, we have $Qu \neq 0$, which implies that $v$ is
proportional to $Qu$.

Let $E$ denote the bundle $\{(H,\ell) \in \sp_{2n} \times \PM^{2n-1}
\mid \im H \subset \ell\} $ over $\PM^{2n-1}$. The first projection
gives a morphism
\[
\pi : E \to \ov\OC_\mini = \OC_\mini \sqcup \{0\}
\]
which is a resolution of singularities. It is an isomorphism over
$\OC_\mini$, and the exceptional fiber is the null section.

The above discussion shows that we have
\[
E \simeq \OC(-1) \otimes_\OC \OC(-1) \simeq \OC(-2).
\]

Let $H^i$ denote $H^i(\PM^{2n-1},\ZM)$, and let
$t \in H^2$ be the first Chern class of $\OC(-1)$.
Then we have $H^*(\PM^{2n-1},\ZM) \simeq \ZM[t]/t^{2n}$.
The Euler class $e$ of $E$ is $2t$.
As $\OC_\mini$ is isomorphic to $E$ minus the null section, we have a
Gysin sequence:
\[
\cdots \longto
H^{i-2} \elem{e} H^i \longto H^i(\OC_\mini,\ZM) \longto
H^{i-1} \elem{e} H^{i+1}
\longto \cdots
\]
As the cohomology of $\PM^{2n-1}$ is concentrated in even degrees,
for $i$ even, we get
$H^i(\OC_\mini,\ZM) \simeq \Coker (e:H^{i-2}\to H^i)$ which is
isomorphic to $\ZM$ for $i = 0$, to $\ZM/2$ if $i$ is an even
integer between $2$ and $4n-2$, and $0$ otherwise. For $i$ odd, we have
$H^i(\OC_\mini,\ZM) \simeq \Ker (e:H^{i-1}\to H^{i+1})$
which is isomorphic to $\ZM$ for $i = 4n - 3$, and zero otherwise.

With $\FM_p$ coefficients, $p$ odd, $\OC_\mini$ has the cohomology of
a sphere. With $\FM_2$ coefficients, it has cohomology
$\FM_2[u]/u^{4n}$, with $u$ in degree one.

Again, $\ov\OC_\mini$ is a cone, so we have
\[
\ic(\ov\OC_\mini,k)_0 = \tau_{\le -1} \circ j_* (\uk_{\OC_\mini}[2n])_0
\simeq \tau_{\le -1} (\RGa(\OC_\mini,k)[2n])
\]
where $j : \OC_\mini \injto \ov\OC_\mini$ is the open immersion.
Thus for $k = \FM_p$, the stalks of $\ic(\ov\OC_\mini,k)$ are as follows: 
\begin{equation*}
\begin{array}{|c|c|c|c|c|c|c|c|}
\hline & -2n & -2n+1 & -2n+2 & \cdots & -2 & -1 & 0 \\
 \hline 
\OC_\mini & k & 0 & 0 & \cdots & 0 & 0 & 0 \\
 \hline 
\{0\} & k & (k)_2 & (k)_2& \cdots & (k)_2 & (k)_2 & 0 \\
 \hline
\end{array}
\end{equation*}

For $k = \ZM$, we get:
\begin{equation*}
\begin{array}{|c|c|c|c|c|c|c|c|}
\hline & -2n & -2n+1 & -2n+2 & \cdots & -2 & -1 & 0 \\
 \hline 
\OC_\mini & \ZM & 0 & 0 & \cdots & 0 & 0 & 0 \\
 \hline 
\{0\} & \ZM & 0 & \ZM/2 & \cdots & \ZM/2 & 0 & 0 \\
 \hline
\end{array}
\end{equation*}
and, for the $p_+$ version, one has to add a copy of $\ZM/2$ in degree $0$ for
the trivial orbit.

Let us give an alternative point of view,
in more concrete terms. If we take $Q = \bigl(
\begin{smallmatrix}
0 & -J\\ J & 0
\end{smallmatrix}
\bigr)$, where
\[
J =
\begin{pmatrix}
0 & & 1\\ & \rdots\\ 1 & & 0
\end{pmatrix}
\]
then we have a morphism $q : \CM^{2n} \to \ov\OC_\mini$ which sends
the vector $(x_1,\dots,x_n,y_1,\dots,y_n)$ to the matrix
\[
\begin{pmatrix}
\begin{matrix}
-x_1y_n & \cdots & -x_1y_1\\ \vdots & \rdots & \vdots\\ -x_ny_n &
\cdots & -x_ny_1\\
\end{matrix}
&
\begin{matrix}
x_1x_n & \cdots & x_1^2\\ \vdots & \rdots & \vdots\\ x_n^2 & \cdots &
x_nx_1
\end{matrix}
\\
\begin{matrix}
-y_1y_n & \cdots & -y_1^2\\ \vdots & \rdots & \vdots\\ y_n^2 & \cdots
& -y_ny_1
\end{matrix}
&
\begin{matrix}
y_1x_n & \cdots & y_1x_1\\ \vdots & \rdots & \vdots\\ y_nx_n & \cdots
& y_nx_1
\end{matrix}
\end{pmatrix}
\]
This can be identified with the quotient by $\{\pm 1\}$. In
particular, this explains why $\ov\OC_\mini$ is $p$-smooth for $p\neq 2$.

The map $q$ induces $q^0 : \CM^{2n}\setminus\{0\} \to \OC_\mini$,
which is the quotient by $\{\pm 1\}$,
which implies that $\OC_\mini\simeq (\CM^{2n}\setminus\{0\})/\{\pm 1\}$
is homotopic to $\RM\PM^{4n-1}$. We have
\[
\RGa(\OC_\mini,k) = \RGa(\RM\PM^{4n-1},k)
= (k \xra{0} k \xra{2} k \xra{0} \dots \xra{2} k \xra{0} k)
\]
where the last complex has a copy of $k$ in each degree between
$0$ and $4n-1$, and the differential is alternatively $0$ and $2$.
We recover the preceding calculation.

Let us note that, for $n = 1$, we recover the ${\sl}_2$ calculation.

\subsection{Subregular class in ${\sl}_n$}
\label{subsec:simple}

Let us consider the nilpotent cone $\NC$ of ${\sl}_n$.
The regular nilpotent orbit $\OC_\reg = \OC_{(n)}$ is open dense in
$\NC$. There is a unique open dense orbit $\OC_\subreg =
\OC_{(n-1,1)}$ in its complement. It is of codimension 2.
Let $U := \OC_\reg \cup \OC_\subreg$.
In this section we will compute the stalks of $\ic(\NC,k)$ restricted 
to $U$ and find a condition on the characteristic of $k$ for
this restriction to be a constant sheaf.

By \cite{BRI,SLO1,SLO2}, the singularity of $\NC$ along $\OC_\subreg$
is a simple surface singularity of type $A_{n-1}$. This is the
singularity at $0$ of the variety $S = \CM^2/\mu_n$, where $\mu_n$ is
the group of $n$th roots of unity in $\CM^*$, and $\zeta\in\mu_n$ acts
on $\CM^2$ by $\left(
\begin{smallmatrix}
\zeta & 0\\
0 & \zeta^{-1}
\end{smallmatrix}
\right)$.
This implies that
$\ic(U,k)_x [-\dim \NC] \simeq \ic(S,k)_0 [-2]$,
where $x\in\OC_\subreg$.
We have
\[
S = \Spec \CM[u,v]^{\mu_n} = \Spec \CM[u^n,v^n,uv] =
\Spec\CM[x,y,z]/(xy-z^n)
\]
and there is a $\CM^*$-action on $\CM^3 = \Spec \CM[x,y,z]$,
contracting to the origin, stabilizing $S$, given by
$t.(x,y,z) = (t^nx,t^ny,t^2z)$. The same argument as with a cone shows
that
\[
(j_* \uk_{S\setminus\{0\}})_0 = \RGa(S\setminus\{0\},k)
= \RGa(S^3/\mu_n,k) = (k \xra{0} k \xra{n} k \xra{0} k)
\]
(a complex in degrees between $0$ and $3$),
where $j:S\setminus\{0\} \injto S$ is the open immersion.

Thus if $k$ is a field of characteristic $p$, then the
stalks of $\ic(U,k)$ are given by:
\[
\begin{array}{|c|c|c|c|c|c|c|}
 \hline 
& - \dim \NC & - \dim \NC + 1 & - \dim \NC + 2 \\
 \hline 
\OC_\reg & k & 0 & 0\\
 \hline
\OC_\subreg & k & (k)_n & 0\\
 \hline
\end{array}
\]
For $k = \ZM$ and for the perversity $p$, we get the constant sheaf on
$X$ in degree $-\dim \NC$. For the perversity $p_+$, we get
an extra stalk $\ZM/n$ in degree $-\dim\NC + 2$ for $\OC_\subreg$.

\begin{remark} Here and in the previous section, 
we have relied on
knowledge of the cohomology ring of real projective spaces and lens 
spaces in order to calculate
the cohomology of spaces like $(\CM^{m} \setminus \{ 0 \}) / \mu_n$.
Thus it is 
not immediately clear how to proceed in the \'etale situation.
However, instead one can use the fact that
$\RGa(\CM^m\setminus\{0\},k)$ is a perfect complex of
$k\mu_n$-modules, because $\mu_n$ acts freely on
$\CM^m\setminus\{0\}$.  As there are only two non-trivial cohomology
groups, there is only one possibility up to quasi-isomorphism. Then
one can take derived invariants to recover
$\RGa((\CM^m\setminus\{0\})/\mu_n,k)$. This proof makes sense if we
replace $\CM$ by $\FM_q$. Alternatively, one could also use comparison
theorems.
\end{remark}

For completeness, and for future use, let us describe the Springer
fiber $\BC_x$. It is the union of $n-1$ projective lines $L_1$, \dots,
$L_{n-1}$, where $L_i$ is identified with the variety of flags
$F_\bullet$ such that $F_j = \im x^{n-1-j}$ for $1 \leq j \leq i-1$
and $F_j = \Ker x^{j-1}$ for $i+1 \leq j \leq n$.  These projective
lines intersect as the Dynkin diagram of type $A_{n-1}$. 
Let $\pi_U : \tilde U \to U$ be the restriction of the Springer
resolution to $U$. Then the stalks of $\pi_{U*}\uk_{\tilde U}[\dim \NC]$ 
are given by:
\begin{equation*}
\begin{array}{|c|c|c|c|c|c|c|}
 \hline 
& - \dim \NC & - \dim \NC + 1 & - \dim \NC + 2 \\
 \hline 
\OC_\reg & k & 0 & 0\\
 \hline
\OC_\subreg & k & 0 & k^{n-1}\\
 \hline
\end{array}
\end{equation*}
When $k$ is of characteristic zero, we can use the decomposition theorem for 
another proof of the above calculation in that case.

Let us note that, for $n = 2$, we recover the ${\sl}_2$ calculation.

\subsection{IC stalks on the nilpotent cone of ${\sl}_3$ for $p \neq 3$}
\label{subsec:sl_3}

Let us consider the nilpotent cone in ${\sl}_3$ and a field of
coefficients $k$ of characteristic $p$ different from $3$. We have
three nilpotent orbits, indexed by the partitions $(3)$, $(21)$ and
$(1^3)$: the regular orbit $\OC_\reg$ of dimension 6, the minimal 
orbit $\OC_\mini$ of dimension 4 and the trivial orbit $\{0\}$.

The complex $\ic(\{0\},k)$ is trivial.
By Subsection \ref{subsec:sl_n}, we know the stalks of
$\ic(\ov\OC_\mini,k)$.
We wish to compute the stalks of $\ic(\NC,k)$.

Let $U := \OC_\reg \cup \OC_\subreg$.
By Subsection \ref{subsec:simple}, we have:
\[
\ic(U, k) = \uk_U [6]
\]
since $k$ is assumed to be of characteristic $p \neq 3$. In other
words, $U$ is $p$-smooth if $p \neq 3$. This
fact will allow us to compute the stalk $\ic(\NC,k)_0$ as in the
preceding cases, where the complement of the origin was smooth. 
Let $j:U\hookrightarrow \NC$ denote the open
immersion. Since $\NC$ is a cone, we have:
\[
(j_*\uk_{U})_0 = \RGa(U,k)
\]

As in Subsection \ref{subsec:simple}, we will consider the
restriction of the Springer resolution $\pi_U : \tilde U \to U$.
We will consider the truncation distinguished triangle:
\[
\tau_{\le -6} \pi_{U*} \uk_{\tilde U}[6] \longto
\pi_{U*}\uk_{\tilde U}[6] \longto
\tau_{> -6} \pi_{U*} \uk_{\tilde U}[6]\rightsquigarrow
\]
By what we said in Subsection \ref{subsec:simple}, 
if $i : \OC_\mini \hookrightarrow U$ denotes
the (closed) inclusion this is isomorphic to the triangle:
\begin{equation}
\label{triangle}
\uk_U[6] \longto \pi_{U*} \uk_{\tilde U}[6] \longto i_* \uk^2_{\OC_\mini}[4] \rightsquigarrow
\end{equation}
Since all the objects involved are actually perverse sheaves,
we have a short exact sequence of perverse sheaves:
\[
0 \longto \uk_U[6] \longto \pi_{U*} \uk_{\tilde U}[6] \longto
i_*\uk^2_{\OC_\mini}[4] \longto 0
\]
with left and right terms intersection cohomology complexes.
We will see that the sequence splits. 
For this, we have to show that
the degree one morphism of the triangle above is zero.
But we have
\[
\begin{array}{rcll}
\Hom(i_*\uk_{\OC_\mini}[4], \uk_U[7])
&=& \Hom(\uk_U[5],i_*\uk_{\OC_\mini}[4])
& \text{by duality}\\
&=& \Hom(i^*\uk_U[5],\uk_{\OC_\mini}[4])
& \text{by adjunction}\\
&=& \Hom(\uk_{\OC_\mini},\uk_{\OC_\mini}[-1])\\
&=& 0.
\end{array}
\]
This is another way of saying that 
$\Ext^1(i_*\uk_{\OC_\mini}[4], \uk_U[6]) = 0$. Thus we have
\[
\pi_{U*} \uk_{\tilde U}[6] = \uk_U[6] \oplus i_*\uk^2_{\OC_\mini}[4].
\]
Taking global sections yields
\[
\RGa(\tilde U,k)[6] = \RGa(U,k)[6] \oplus \RGa(\OC_\mini,k)^2[4].
\]

Finally if we can compute the cohomology of $\tilde U$ and $\OC_\mini$ with
coefficients in $k$, then we will have computed the cohomology of $U$
and equivalently, as $U$ is $p$-smooth, its intersection cohomology.

Recall that $\tilde U$ is the cotangent bundle to the flag variety $G/B$ with
the zero section removed.  We can then apply the Gysin sequence
\[
\cdots \longto H^{i-6}\elem{e} H^i \longto H^i(\tilde U,k) \longto
H^{i-5} \elem{e} H^{i+1} \longto \cdots
\]
where $H^i := H^i(G/B,k)$
and $e$ denotes multiplication by the Euler class of $T^*(G/B)$, which in this
case is given by $\chi(G/B) = 6$ times a generator of $H^6$.  Recalling that
the cohomology of $G/B$ is as follows:

\begin{center}
\begin{tabular}{ccccccc}
0 & 1 & 2 & 3 & 4 & 5 & 6 \\ \hline $k$ & 0 & $k^2$ & 0 & $k^2$ & 0 &
$k$ 
\end{tabular}
\end{center}
we find the cohomology of of the complex $\RGa(\tilde U,k)[6]$ is given by:
\begin{center}
\begin{tabular}{ccccccccccccc}
-6 & -5 & -4 & -3 & -2 & -1 & 0 & 1 & 2 & 3 & 4 & 5 & 6 \\ \hline $k$
& 0 & $k^2$ & 0 & $k^2$ & $(k)_6$ & $(k)_6$ & $k^2$ & 0 & $k^2$ & 0 &
$k$ & 0
\end{tabular}
\end{center}

We have computed the cohomology of $\OC_\mini$ in Subsection
\ref{subsec:sl_n}. Since $p \neq 3$, the cohomology of
the complex $\RGa(\OC_\mini,k)[4]$ is given by:
\[
\begin{array}{cccccccc}
-4 & -3 & -2 & -1 & 0 & 1 & 2 & 3 \\
\hline
k & 0 & k & 0 & 0 & k & 0 & k
\end{array}
\]

Subtracting two copies of the cohomology of $\RGa(\OC_\mini,k)[4]$ from the
cohomology of $\RGa(\tilde U,k)[6]$, we obtain the cohomology of 
$\RGa(U,k)[6]$ for $p\neq 3$:

\begin{center}
\begin{tabular}{ccccccccccccc}
-6 & -5 & -4 & -3 & -2 & -1 & 0 & 1 & 2 & 3 & 4 & 5 & 6 \\ \hline $k$
& 0 & 0 & 0 & 0 & $(k)_2$ & $(k)_2$ & 0 & 0 & 0 & 0 & $k$ & 0
\end{tabular}
\end{center}
By truncation, we obtain the stalk $\ic(\NC,k)_0$ for $p\neq 3$:
\begin{thm} Let $k$ be a field of characteristic $p \ne 3$. Then
\[
\ic(\NC,k)_0 = 
\begin{cases}
k[6]\oplus k[1] & \text{if $p = 2$,} \\
k[6] & \text{if $p > 3$.}
\end{cases}
\]
\end{thm}

This completes the geometric determination of the decomposition matrix of
$GL_3$-equivariant perverse sheaves on its nilpotent cone, and hence
of the symmetric group $\SG_3$, in characteristic $2$: 
\[
\begin{array}{cc}
& 1^3 \ 21 \ 3\\
\begin{array}{c}
1^3\\
21\\
3
\end{array}
&
\begin{pmatrix}
{\bf 1}&{\bf 0}&0\\
{\bf 0}&{\bf 1}&0\\
{\bf 1}&{\bf 0}&1
\end{pmatrix}
\end{array}
\qquad
\begin{array}{cc}
& D_{(3)} \ D_{(21)} \\
\begin{array}{c}
S_{(3)}\\
S_{(21)}\\
S_{(1^3)}
\end{array}
&
\begin{pmatrix}
&{\bf 1}&{\bf 0}&\\
&{\bf 0}&{\bf 1}&\\
&{\bf 1}&{\bf 0}&
\end{pmatrix}
\end{array}
\]

\subsection{Some IC stalks on the nilpotent cone of ${\sl}_4$ for $p\neq 2$}

In this subsection, we calculate the stalks of the 
intersection cohomology complex
on $\ov\OC_{(2^2)}$ in characteristic $p \ne 2$. The
closure $\ov\OC_{(2^2)}$ consists of three orbits: $\OC_{(2^2)}$, 
$\OC_\mini = \OC_{(2, 1^2)}$ and $\{ 0 \}$ of dimensions 
8, 6 and 0. Let $U' = \OC_{(2^2)} \cup \OC_\mini$.

Kraft and Procesi show~\cite{kp} that the singularity of
$U'$ along $\OC_\mini$ is equivalent to a simple $A_1$ singularity,
that of the nilpotent cone of ${\sl_2}$ that we studied in
Section \ref{subsec:sl_2}. It follows that the
stalks of $\ic(\ov\OC_{(2^2)},k)$ restricted to $U'$ are as follows:

\begin{equation*}
\begin{array}{|c|c|c|c|c|c|c|}
 \hline 
& - 8 & - 7 & - 6 \\
 \hline 
\OC_{(2^2)} & k & 0 & 0\\
 \hline
\OC_\mini & k & (k)_2 & 0\\
 \hline
\end{array}
\end{equation*}

>From now on, let $k$ be a field of characteristic $p \ne 2$ (so that $U'$ is $p$-smooth).
Using the Deligne construction and the fact that $\ov\OC_{(2^2)}$
is a cone we have
\[
\ic(\ov\OC_{(2^2)}, k)_0 = ( \tau_{\le -1} j_* \uk_{U'}[8])_0 
= \tau_{\le -1}(\RGa(U',k)[8])
\]
where $j : U' \hookrightarrow \ov\OC_{(2^2)}$ denotes the inclusion.
We will proceed as in the previous section, calculating 
$\RGa(U',k)$ with the help of the resolution
\[
\pi : T^* (G/P) \to \ov\OC_{(2^2)}
\]
defined in Section~\ref{subsec:richardson}. Here $G = SL_4$ and $P$ 
is the parabolic stabilizing $\CM^2 \subset \CM^4$. The fiber over a
nilpotent $N$ with Jordan blocks $(2,1^2)$ is the collection of
2-planes $V\in \Gras(2,4)$ which are in the kernel of $N$ and contain its
image.  The image of $N$ is 1-dimensional and contains the kernel which is
3-dimensional.  Thus the collection of such 2-planes forms a
projective line $\PM^1$. Similarly to above, let $\tilde{U'}$ denote
the preimage of $U'$ under $\pi$ and denote by $\pi_{U}$ the 
restriction of $\pi$ to $U'$.

We consider the push-forward
$\pi_{U*}\uk_{U'}[8]$. The last paragraph implies that the
stalks of this push-forward are as follows:
\begin{equation*}
\begin{array}{|c|c|c|c|c|c|c|}
 \hline 
& - 8 & - 7 & - 6 \\
 \hline 
\OC_\reg & k & 0 & 0\\
 \hline
\OC_\subreg & k & 0 & k\\
 \hline
\end{array}
\end{equation*}
Let $i:\OC_\mini\hookrightarrow U$ denote the closed inclusion. We 
have a truncation distinguished triangle analogous to 
Section \ref{subsec:sl_3}:
\[
\uk_{U'}[8] \longto \pi_{U*} \uk_{\tilde {U'}}[8] \longto i_* \uk_{\OC_\mini}[6] \rightsquigarrow
\]
or equivalently a short exact sequence of perverse sheaves:
\[
0 \longto \uk_{U'}[8] \longto \pi_{U*} \uk_{\tilde {U'}}[8] \longto
i_*\uk_{\OC_\mini}[6] \longto 0
\]
This sequence splits by the same calculation as in Section \ref{subsec:sl_3}.
Hence
\[
\pi_{U*} \uk_{\tilde {U'}}[8] = \uk_{U'}[8] \oplus i_*\uk_{\OC_\mini}[6]
\]
and
\begin{equation} \label{eq:splitsl4}
\RGa(\tilde {U'},k)[8] = \RGa(U',k)[8] \oplus \RGa(\OC_\mini,k)[6].
\end{equation}

We already know the cohomology of $\OC_{\mini}$ by Section 
\ref{subsec:sl_n}. Hence all that remains is to compute $\RGa(\tilde {U'},k)[8]$.
Recall that the cohomology of $\Gras(2,4)$ has a integral basis given by
Schubert cells. In particular its cohomology over $\ZM$ is:
\begin{center}
\begin{tabular}{ccccccccc}
0 & 1 & 2 & 3 & 4 & 5 & 6 & 7 & 8 \\ \hline 
$\ZM$ & 0 & $\ZM$ & 0 & $\ZM^2$ & 0 & $\ZM$ & 0 & $\ZM$
\end{tabular}
\end{center}
 The space $\tilde{U'}$ is the cotangent bundle to $\Gras(2,4)$ with the 
zero section removed.  Using the Gysin
sequence we can compute the cohomology of $\tilde U'$ from that of
$\Gras(2,4)$. Thus the cohomology of $\uk_{\tilde{U'}}[8]$ is given by:

\begin{center}
\begin{tabular}{cccccccccccccccc}
-8 & -7 & -6 & -5 & -4 & -3 & -2 & -1 & 0 & 1 & 2 & 3 & 4 & 5 & 6 &
7\\ \hline $k$ & 0 & $k$ & 0 & $k^2$ & 0 & $k$ & $(k)_6$ & $(k)_6$ &
$k$ & 0 & $k^2$ & 0 & $k$ & 0 & $k$
\end{tabular}
\end{center}

We computed the cohomology of $\OC_\mini$ in Section~\ref{subsec:sl_n}
(shifted here by 6) for $p\neq 2$ to be:
\begin{center}
\begin{tabular}{cccccccccccc}
-6 & -5 & -4 & -3 & -2 & -1 & 0 & 1 & 2 & 3 & 4 & 5 \\ \hline $k$ & 0
& $k$ & 0 & $k$ & 0 & 0 & $k$ & 0 & $k$ & 0 & $k$
\end{tabular}
\end{center}

By (\ref{eq:splitsl4}) we obtain the desired cohomology of $U^{\prime}$ 
(for $p\neq 2$) shifted by 8:

\begin{center}
\begin{tabular}{cccccccccccccccc}
-8 & -7 & -6 & -5 & -4 & -3 & -2 & -1 & 0 & 1 & 2 & 3 & 4 & 5 & 6 &
7\\ \hline $k$ & 0 & 0 & 0 & $k$ & 0 & 0 & $(k)_3$ & $(k)_3$ & 0 & 0 &
$k$ & 0 & 0 & 0 & $k$
\end{tabular}
\end{center}

Finally, all that remains is to truncate:

\begin{thm} Let $k$ be a field of characteristic $p \ne 2$. Then
\[
\ic(\ov\OC_{(2^2)},k)_0 = 
\begin{cases}
k[8]\oplus k[4] \oplus k[1] & \text{if $p = 3$,} \\
k[8] \oplus k[4] & \text{if $p > 3$.}
\end{cases}
\]
\end{thm}

We get the following parts of the decomposition matrices for $p = 3$:
\[
\begin{array}{cc}
& 1^4 \ 21^2 \ 2^2\\
\begin{array}{c}
1^4\\
21^2\\
2^2
\end{array}
&
\begin{pmatrix}
{\bf 1}&{\bf 0}&{\bf 0}\\
{\bf 0}&{\bf 1}&{\bf 0}\\
{\bf 1}&{\bf 0}&{\bf 1}
\end{pmatrix}
\end{array}
\qquad
\begin{array}{cc}
& D_{(4)} \ D_{(31)} \ D_{(2^2)}\\
\begin{array}{c}
S_{(4)}\\
S_{(31)}\\
S_{(2^2)}
\end{array}
&
\begin{pmatrix}
&{\bf 1}&{\bf 0}&{\bf 0}&\\
&{\bf 0}&{\bf 1}&{\bf 0}&\\
&{\bf 1}&{\bf 0}&{\bf 1}&
\end{pmatrix}
\end{array}
\]

\newcommand{\etalchar}[1]{$^{#1}$}

\end{document}